\documentclass[a4paper, reqno]{amsart}
\usepackage[super,compress]{cite}
\usepackage{latexsym}
\usepackage{amsmath}
\usepackage{amsfonts}
\usepackage{amssymb}
\usepackage{amsxtra}
\usepackage{graphicx,comment,color,mathrsfs}
\usepackage{breakurl}
\usepackage[all]{xy}

\newtheorem{definition}{Definition}[section]

\newtheorem{theorem}[definition]{Theorem}

\newtheorem{remark}[definition]{Remark}

\newcommand{\Z}{\mathbb{Z}}
\newcommand{\C}{\mathbb{C}}

\newcommand{\Com}{\operatorname{Com}}

\newcommand{\Q}{\mathbb{Q}}
\newcommand{\g}{\mathfrak{g}}
\newcommand{\h}{\mathfrak{h}}
\newcommand{\ns}{\mathfrak{ns}}
\newcommand{\sll}{\mathfrak{sl}}
\newcommand{\gl}{\mathfrak{gl}}

\newcommand{\poch}[1]{\left(#1;q\right)_{\infty}}

\newcommand{\ssqrt}[1]{\operatorname{\sqrt{\smash[b]{#1}}}}

\newcommand{\W}{\mathcal{W}}

\begin{document}
\markboth{Shigenori Nakatsuka}{Feigin--Semikhatov conjecture}

\title{Feigin--Semikhatov conjecture and related topics
}

\author{Shigenori Nakatsuka
}

\address{Department of Mathematical and Statistical Sciences, University of Alberta, CAB 632, Edmonton, AB T6G 2G1, Canada}
\email{shigenori.nakatsuka@ualberta.ca}

\maketitle

\begin{abstract}

Feigin--Semikhatov conjecture, now established, states algebraic isomorphisms between the cosets of the subregular $\W$-algebras and the principal $\W$-superalgebras of type A by their full Heisenberg subalgebras. It can be seen as a variant of Feigin--Frenkel duality between the $\W_n$-algebras and also as a generalization of the connection between the $\mathcal{N}=2$ superconformal algebra and the affine algebra $\widehat{\sll}_{2,k}$. 

We review the recent developments on the correspondence of the subregular $\W$-algebras and the principal $\W$-superalgebras of type A at the level of algebras, modules and intertwining operators, including fusion rules.

\end{abstract}

\section{Introduction}	

Vertex (super)algebras are symmetry algebras (or chiral algebras) for two dimensional conformal field theories and have been introduced in the mathematical literature by Borcherds\cite{B}. 
They also appear as algebras of local operators in several higher dimensional quantum field theories, which have attracted an intensive attention in the last decade, see Refs. \citen{AGT,BRcR,GR,CG,FG,FeiG,CDGG,CCFFGHP} for example.

The $\W$-superalgebras and their affine cosets, that is, subalgebras whose fields commute with affine subalgebras inside give a rich class of such vertex algebras. The $\W$-superalgebras, denoted by $\W^k(\g,f)$, are obtained from the affine vertex superalgebra $V^k(\g)$ through the quantum Drinfeld--Sokolov reduction parametrized by nilpotent elements $f$ inside $\g$. In this paper, we will write $\W_k(\g,f)$ and $L_k(\g)$ for their simple quotients. When $\g$ is a simply-laced Lie algebra and $f=f_{\mathrm{prin}}$ is called principal (or regular), $\W^k(\g)=\W^k(\g,f_{\mathrm{prin}})$ enjoys the triality consisting of the Feigin--Frenkel duality \cite{FF} and the Goddard--Kent--Olive type construction \cite{GKO,ACL}--- an example of affine cosets---, which asserts that all of them give the same algebra. 
In particular when $\g=\sll_n$, the triality for $\W^k(\sll_n)$, also known as Fateev--Lukyanov's $\W_n$-algebra \cite{FL}, is depicted in Fig. \ref{Triality for principal} below.
\begin{figure}[htbp]
\centering
\caption{}
\setlength{\unitlength}{1mm}
\begin{picture}(90,17)(2,0)\label{Triality for principal} 
\put(2,10){\footnotesize$\W^{k}(\sll_n)$}
\put(47,10){\footnotesize$\W^{\check{k}}(\sll_n)$}
\put(8,0){\footnotesize$\Com(V^\ell(\sll_n),V^{\ell-1}(\sll_n)\otimes L_1(\sll_n))$}
\put(15,11){\line(1,0){29}}
\put(15,7){\line(3,-1){10}}
\put(45,7){\line(-3,-1){10}}
\put(24,12){\footnotesize{FF duality}}
\put(20,6){\footnotesize{GKO}}
\put(80,0){\footnotesize $\tfrac{1}{k+n}+\tfrac{1}{\ell+n}=1$}
\put(80,5){\footnotesize $(k+n)(\check{k}+n)=1$}
\put(80,10){\footnotesize \underline{Relation of levels}}
\end{picture}
\end{figure}

Another kind of this tirality was suggested first by Feigin and Semikhatov \cite{FS} for the $\W$-algebra $\W^k(\sll_n,f_{\mathrm{sub}})$ associated with the subregular nilpotent element $f_{\mathrm{sub}}$ as in Fig. \ref{Triality for subregular}.

\begin{figure}[htbp]
\centering
\caption{}
\setlength{\unitlength}{1mm}
\begin{picture}(90,17)(2,0)\label{Triality for subregular} 
\put(-9,10){\footnotesize$\Com(\pi,\W^{k}(\sll_n,f_{\mathrm{sub}}))$}
\put(41,10){\footnotesize$\Com(\pi,\W^{\check{k}}(\sll_{n|1},f_{\mathrm{prin}}))$}
\put(14,0){\footnotesize$\Com(V^\ell(\gl_n),V^{\ell}(\sll_{n|1}))$}
\put(23,11){\line(1,0){16}}
\put(15,7){\line(3,-1){10}}
\put(45,7){\line(-3,-1){10}}
\put(80,0){\footnotesize $\tfrac{1}{k+n}+\tfrac{1}{\ell+n}=1$}
\put(80,5){\footnotesize $(k+n)(\check{k}+n-1)=1$}
\put(80,10){\footnotesize \underline{Relation of levels}}
\end{picture}
\end{figure}
Here we take the coset of $\W$-superalgebras by the maximal affine subalgebra inside, which is rank one Heisenberg algebra $\pi$. 
The very beginning case of the Feigin--Frenkel type duality in this family is indeed derived from the so-called Kazam--Suzuki coset construction \cite{KaSu} of the $\mathcal{N}=2$ superconformal algebra.

Gaiotto and Rap\v{c}\'{a}k\cite{GR} found a vast generalization of these trialities for \emph{vertex algebras at the corner} appearing as local operators at the junctions of supersymmetric interfaces in $\mathcal{N} = 4$ Super Yang--Mills gauge theory. 
The above examples correspond to the following diagrams which describe the ranks of the gauge group $U(n)$ and the half-BPS interfaces placed as three dimensional boundary conditions.
\setlength{\unitlength}{1mm}
\begin{figure}[htbp]
\centering
\caption{}
\begin{picture}(90,35)(0,0)
\put(0,23){\line(1,0){6}}
\put(0,23){\line(0,1){6}}
\put(0,23){\line(-1,-1){4}}
\put(2,25){\footnotesize$n$}
\put(26,23){\line(1,0){6}}
\put(26,23){\line(0,1){6}}
\put(26,23){\line(-1,-1){4}}
\put(27,19){\footnotesize$n$}
\put(14,7){\line(1,0){6}}
\put(14,7){\line(0,1){6}}
\put(14,7){\line(-1,-1){4}}
\put(9,8){\footnotesize$n$}
\put(12,23){\vector(1,0){4}}
\put(25,16){\vector(-3,-1){4}}
\put(6,15){\vector(-3,1){4}}
\put(60,23){\line(1,0){6}}
\put(60,23){\line(0,1){6}}
\put(60,23){\line(-1,-1){4}}
\put(62,25){\footnotesize$n$}
\put(61,19){\footnotesize$1$}
\put(86,23){\line(1,0){6}}
\put(86,23){\line(0,1){6}}
\put(86,23){\line(-1,-1){4}}
\put(87,19){\footnotesize$n$}
\put(82,24){\footnotesize$1$}
\put(74,7){\line(1,0){6}}
\put(74,7){\line(0,1){6}}
\put(74,7){\line(-1,-1){4}}
\put(69,8){\footnotesize$n$}
\put(76,9){\footnotesize$1$}
\put(72,23){\vector(1,0){4}}
\put(85,16){\vector(-3,-1){4}}
\put(66,15){\vector(-3,1){4}}
\end{picture}
\end{figure}
\vspace{-6mm}
A large class of this triality is now proven by Creutzig and Linshaw \cite{CL1,CL2} mathematically and its application to the representation theory of $\W$-superalgebras is also under work in progress \cite{CGN,CGNS,CLNS}.

In this article, we restrict ourselves to the case considered by Feigin and Semikhatov and review the results obtained at this time with an emphasis on how the Kazama--Suzuki coset construction is generalized and refined in this setting.
\section*{Acknowledgments}
The author would like to express his deepest gratitude to Thomas Creutzig for collaboration and valuable comments during his preparation of this article. He also thanks Naoki Genra, Andrew. R. Linshaw and Ryo Sato for collaboration and for Masahito Yamazaki and Yuto Moriwaki for useful discussions.
The author is supported by JSPS  KAKENHI Grant Number 20J1014 and JSPS Overseas Research Fellowships Grant Number 202260077.
The work is supported by World Premier International Research Center Initiative (WPI Initiative), MEXT, Japan.

\section{The case of $\widehat{\sll}_{2,k}$ and $\mathfrak{ns}_{2,c}$}
\subsection{Coset construction of $\mathfrak{ns}_{2,c}$}

The $\mathcal{N}=2$ superconformal algebra $\mathfrak{ns}_{2,c}$ with central charge $c$ has four generating fields: $L(z)$, $J(z)$ of even parity and $G^\pm(z)$ of odd parity satisfying the following operator product expansions (OPEs)
\begin{align*}
&L(z)L(w)\sim \frac{\tfrac{1}{2}c}{(z-w)^4}+\frac{2L(w)}{(z-w)^2}+\frac{\partial L(w)}{z-w},\quad  J(z)J(w)\sim \frac{\tfrac{1}{3}c}{(z-w)^2},\\
&L(z)J(w)\sim \frac{J(w)}{(z-w)^2}+\frac{\partial J(w)}{(z-w)},\quad L(z)G^\pm(w)\sim \frac{\frac{3}{2}G^\pm(w)}{z-w}+\frac{\partial_wG^\pm(w)}{z-w},\\
&J(z)G^\pm(w) \sim\frac{\pm G^\pm(w)}{(z-w)},\quad  G^\pm(z)G^\pm(w)\sim0,\\
&G^\pm(z)G^\mp(w)\sim \frac{\tfrac{2}{3}c}{(z-w)^3}+\frac{2J(w)}{(z-w)^2}+\frac{2L(w)+\partial J(w)}{z-w}.
\end{align*}
For clarification, let $V^c(\ns_2)$ denote the (universal) vertex algebra freely generated by these fields and by $L_c(\ns_2)$ its unique simple quotient. 
Kazama and Suzuki\cite{KaSu} gave a family of realizations of $V^c(\ns_2)$ based on Hermitian symmetric spaces $G/H$.
We are interested in the case $G/H=SU(n+1)/SU(n)\times U(1)$. 
More precisely, we consider the vertex superalgebra $V^k(\sll_{n+1})\otimes bc^{\otimes n}$ of the universal affine vertex algebra of $\sll_{n+1}$ at level $k$ and $n$ copies of $bc$-systems. 
Since we have homomorphisms $V^k(\gl_n) \rightarrow V^k(\sll_{n+1})$ and $V^1(\gl_n)\rightarrow bc^{\otimes n}$, we may take the diagonal coset $\Com(V^{k+1}(\gl_n), V^k(\sll_{n+1})\otimes bc^{\otimes n})$, i.e., the subalgebra consisting of elements whose fields commute with the diagonal $\widehat{\gl}_{n,k+1}$-action on $V^k(\sll_{n+1})\otimes bc^{\otimes n}$. 
Then the coset contains $V^c(\ns_2)$ with specific central charge $c$ whereas the whole coset is conjectured\cite{I,CL0} to be isomorphic to the principal $\W$-superalgebra $\W(\sll_{n+1|n})$ at certain level:
\begin{align}
V^c(\ns_2)\hookrightarrow \W(\sll_{n+1|n})\simeq \Com\left(V^{k+1}(\gl_n),V^k(\sll_{n+1})\otimes bc^{\otimes n}\right).
\end{align}
When $n=1$, $V^c(\ns_2)$ is exactly the principal $\W$-superalgebra $\W(\sll_{2|1})$ and gives the honest coset construction of $V^c(\ns_2)$.
More explicitly, the isomorphism is given by the following:
\begin{align}\label{eq KS}
\begin{array}{cccc}
{\bf KS}\colon&V^c(\ns_{2})&\xrightarrow{\simeq}&\Com\left(\pi^{H_+^{\Delta}},V^k(\mathfrak{sl}_2) \otimes V_{\Z}\right)\\
&G^+(z)&\mapsto& \sqrt{\frac{2}{k+2}}e(z)\otimes V_{1}(z)\\
&G^-(z)&\mapsto& \sqrt{\frac{2}{k+2}}f(z)\otimes V_{-1}(z)\\
&J(z)&\mapsto &\frac{-2}{k+2}H_+^{\Delta}(z)+1 \otimes b_1(z)
\end{array}
\end{align}
with 
\begin{align}\label{central charge}
H_+^{\Delta}(z)=-\frac{1}{2}h(z)\otimes 1+1\otimes b_1(z),\quad c=\frac{3k}{k+2}.
\end{align}
Here we have used the boson-fermion correspondence between the $bc$-system and the lattice vertex superalgebra $V_\Z$ associated with the lattice $\Z$; $V_n(z)$ is the vertex operator for the element in the lattice $n\in \Z$; $b_n(z)$ is the Heisenberg field satisfying the OPE $b_n(z)b_n(z)\sim n^2/(z-w)^2$; $\pi^{H_+^\Delta}$ is the Heisenberg vertex algebra generated by $H_+^{\Delta}(z)$. 
The construction \eqref{eq KS} descends to the simple quotients
\begin{align}
L_c(\ns_2)\xrightarrow{\simeq}\Com\left(\pi^{H_+^{\Delta}},L_k(\sll_2) \otimes V_{\Z}\right).
\end{align}
When $k$ is a non-negative integer, this realization has been used to construct the unitary minimal representations \cite{DVPYZ} as a variant of the celebrated Goddard--Kent--Olive construction\cite{GKO} of those representations for the Virasoro algebra. We will return back to this point later.

The coset construction \eqref{eq KS} implies that the representation theories of $V^c(\ns_2)$ and $V^k(\sll_2)$ share some similarity. This means that one can not only go from the $V^k(\sll_2)$-side to the $V^c(\ns_2)$-side, but also in the opposite direction. 
Feigin--Semikhatov--Tipunin \cite{FST} found the crucial coset type realization of $V^k(\sll_2)$ in terms of $V^c(\ns_2)$ by using the negative-definite lattice vertex superalgebra $V_{\ssqrt{-1}\Z}$. Namely, they established the following isomorphism of vertex  algebras:
\begin{align}\label{eq FST}
\begin{array}{cccc}
{\bf FST}\colon& V^k(\mathfrak{sl}_2)&\xrightarrow{\simeq}&\Com(\pi^{H_-^{\Delta}},V^c(\ns_{2}) \otimes V_{\ssqrt{-1}\Z})\\
&e(z)&\mapsto&\sqrt{\frac{k+2}{2}}G^+(z)\otimes V_{\ssqrt{-1}}(z)\\
&f(z)&\mapsto&\sqrt{\frac{k+2}{2}}G^-(z)\otimes V_{\text{-}\ssqrt{-1}}(z)\\
&\frac{1}{2}h(z)&\mapsto &\frac{k+2}{2}H_-^{\Delta}(z)-1\otimes b_{\ssqrt{-1}}(z)
\end{array}
\end{align}
with $H_-^\Delta(z)=J(z)\otimes 1+1\otimes b_{\ssqrt{-1}}(z)$. Again, \eqref{eq FST} descends to the simple quotients
\begin{align}
L_k(\sll_2)\xrightarrow{\simeq}\Com\left(\pi^{H_-^{\Delta}},L_c(\ns_{2}) \otimes  V_{\ssqrt{-1}\Z}\right).
\end{align}

\subsection{Equivalence of module categories}
Let us recall the block-wise equivalence \cite{FST,Sato} of representation categories for $V^k(\sll_2)$ and $V^c(\ns_2)$ by the coset constructions \eqref{eq KS} and \eqref{eq FST}. 
From now on, we always assume $c=\frac{3k}{k+2}$ as in \eqref{central charge}. 
Since the coset constructions themselves use Heisenberg fields, it is reasonable to consider those modules of  $V^k(\sll_2)$ and $V^c(\ns_2)$ which behave nicely for their Heisenberg fields inside them. This motivates to introduce the category of weight modules: let $V^k(\sll_2)\text{-mod}_{\mathrm{wt}}$ denote the category of $V^k(\sll_2)$-modules which decompose into direct sums of Fock modules with respect to $\pi^{\frac{1}{2}h}$. Note that the eigenvalues of $\frac{1}{2}h_{0}$ on $V^k(\sll_2)$ is $\Z$. Then the category $V^k(\sll_2)\text{-mod}_{\mathrm{wt}}$ admits a block decomposition in terms of the spectrum of $\frac{1}{2}h_{0}$ modulo $\Z$:
\begin{align*}
V^k(\sll_2)\text{-mod}_{\mathrm{wt}}=\bigoplus_{[\lambda]\in \C/\Z}V^k(\sll_2)\text{-mod}_{\mathrm{wt}}^{[\lambda]}.
\end{align*}
Similarly, we introduce the category $V^c(\ns_2)\text{-mod}_{\mathrm{wt}}$ of weight $V^c(\ns_2)$-modules in terms of $\pi^J$. Then it decomposes into 
\begin{align*}
V^c(\ns_2)\text{-mod}_{\mathrm{wt}}=\bigoplus_{[\mu]\in \C/\Z}V^c(\ns_2)\text{-mod}_{\mathrm{wt}}^{[\mu]}
\end{align*}
by using the spectrum of $J_{0}$ on $V^c(\ns_2)$, which is $\Z$. To compare the categories of weight modules, it suffices to notice that the coset constructions for algebras \eqref{eq KS} and \eqref{eq FST} are generalized to modules by using the multiplicity spaces of suitable Fock modules. To obtain a functor from $V^k(\sll_2)\text{-mod}_{\mathrm{wt}}$ to $V^c(\ns_2)\text{-mod}_{\mathrm{wt}}$, take an arbitrary weight module $M$ and decompose $M\otimes V_\Z$ as a module over $V^c(\ns_2)\otimes \pi^{H_+^{\Delta}}$:
\begin{align}
M\otimes V_\Z\simeq \bigoplus_{\xi\in \C}\Omega^+_{\xi}(M)\otimes \pi^{H_+^{\Delta}}_\xi
\end{align}
with 
\begin{align*}
\Omega_{\xi}^+(M):=\{a\in M\otimes V_\Z\mid H_{+,n}^{\Delta}a=\delta_{n,0}\xi \ a,\ (n\geq0)\}
\end{align*}
and $\pi^{H_+^{\Delta}}_\xi$ the Fock module of $\pi^{H_+^{\Delta}}$ on which $H_{+,0}^{\Delta}$ acts by $\xi$. If $M$ lies in $V^k(\sll_2)\text{-mod}_{\mathrm{wt}}^{[\lambda]}$, then $\Omega_{\xi}^+(M)$ lies in $V^c(\ns_2)\text{-mod}_{\mathrm{wt}}^{[-\varepsilon \xi]}$ with $\varepsilon=\frac{2}{k+2}$. Therefore, $\Omega_\xi^+$ defines a functor 
\begin{align}
\Omega_\xi^+\colon V^k(\sll_2)\text{-mod}_{\mathrm{wt}}^{[\lambda]}\rightarrow V^c(\ns_2)\text{-mod}_{\mathrm{wt}}^{[-\varepsilon \xi]}.
\end{align}
It is non-zero if and only if $\xi\in -\lambda+\Z$.
To get a functor in the opposite direction, take an arbitrary $V^c(\ns_2)$-module $N$ and decompose $N\otimes V_{\ssqrt{-1}\Z}$ as a module over $V^k(\sll_2)\otimes \pi^{H_-^{\Delta}}$:
\begin{align}
N\otimes V_{\ssqrt{-1}\Z}\simeq \bigoplus_{\xi\in \C}\Omega^-_{\xi}(N)\otimes \pi^{H_-^{\Delta}}_\xi.
\end{align}
Then we obtain a functor 
\begin{align}
\Omega^-_{\xi}\colon V^c(\ns_2)\text{-mod}_{\mathrm{wt}}^{[\mu]}\rightarrow V^k(\sll_2)\text{-mod}_{\mathrm{wt}}^{[\varepsilon^{-1}\xi]},
\end{align}
which is non-zero if and only if $\xi\in \lambda+\Z$.

The functors $\Omega^\pm_\xi$ are either zero or an equivalence. In particular, the block-wise equivalence of categories is stated in the following way:
\begin{theorem}\cite{FST,Sato}\label{equivalence of categories}
The functors 
\begin{align}
\Omega_{-\lambda}^+\colon V^k(\sll_2)\text{-}\mathrm{ mod}_{\mathrm{wt}}^{[\lambda]}\ \rightleftarrows  \ V^c(\ns_2)\text{-}\mathrm{mod}_{\mathrm{wt}}^{[\varepsilon \lambda]}\colon \Omega^-_{\varepsilon\lambda}
\end{align}
are quasi-inverse to each other and thus give an equivalence of categories. 
\end{theorem}
The proof is established by giving natural isomorphisms connecting modules $M$ and their images in $M\otimes V_\Z\otimes V_{\ssqrt{-1}\Z}$ under the compositions $\Omega_{-\lambda}^+ \circ \Omega_{\varepsilon \lambda}^-$ and $\Omega_{\varepsilon \lambda}^-\circ \Omega_{-\lambda}^+$.

\subsection{Some basic modules and resolutions}\label{Some basic modules and resolutions}
To illustrate the block-wise equivalence, we compare resolutions of simple modules by some basic modules when the level $k$ is irrational. On $V^k(\sll_2)$-side, the Weyl module $\mathbb{V}^k_{n\varpi_1}$ of highest weight $n\varpi_1$, the induced module whose top space is the $n+1$-dimensional $\sll_2$-module, has a two step resolution by the affine Verma modules:

\begin{align}\label{BGG type resolution}
0\rightarrow \mathbb{M}^k_{-(n+2)\varpi_1}\rightarrow \mathbb{M}^k_{n\varpi_1}\rightarrow \mathbb{V}^k_{n\varpi_1}\rightarrow 0.
\end{align}
Applying the functor $\Omega_{-n/2}^+$, we obtain the following resolution of the simple $V^c(\ns_2)$-module $L_c(\frac{1}{4}\varepsilon n,\frac{1}{2}\varepsilon n)$ of highest weight $(L_0,J_0)=(\frac{1}{4}\varepsilon n,\frac{1}{2}\varepsilon n)$:
\begin{align*}
0\rightarrow \Omega^+_{-n/2}(\mathbb{M}^k_{-(n+2)\varpi_1})\rightarrow \Omega^+_{-n/2}(\mathbb{M}^k_{n\varpi_1})\rightarrow L_c(\tfrac{1}{4}\varepsilon n,\tfrac{1}{2}\varepsilon n)\rightarrow0.
\end{align*}
Each component has the following character $\mathrm{tr}_{\bullet}(q^{L_0}z^{J_0})$:
\begin{align*}
&\mathrm{ch}\ L_c(\tfrac{1}{4}\varepsilon n,\tfrac{1}{2}\varepsilon n)=q^{\frac{1}{4}\varepsilon n}z^{\frac{1}{2}\varepsilon n}\frac{1-q^{n+1}}{1+z^{-1}q^{n+\frac{1}{2}}}
\frac{\poch{-zq^{\frac{3}{2}},-z^{-1}q^{\frac{1}{2}} }}{\poch{q}^2},\\
&\mathrm{ch}\ \Omega^+_{-n/2}(\mathbb{M}^k_{n\varpi_1})=q^{\frac{1}{4}\varepsilon n}z^{\frac{1}{2}\varepsilon n}
\frac{\poch{-zq^{\frac{3}{2}},-z^{-1}q^{\frac{1}{2}} }}{\poch{q}^2},\\
&\mathrm{ch}\ \Omega^+_{-n/2}(\mathbb{M}^k_{-(n+2)\varpi_1})=q^{\frac{1}{4}\varepsilon n+\frac{1}{2}(n+1)^2}z^{\frac{1}{2}\varepsilon n-(n+1)}
\frac{\poch{-zq^{\frac{1}{2}-n},-z^{-1}q^{\frac{3}{2}+n} }}{\poch{q}^2},
\end{align*}
in terms of the $q$-Pochhammer symbols
$(a_1,\cdots,a_m;q)_\infty=\prod_{i=1}^m\prod_{N=0}^\infty(1-a_iq^{N})$.
As indicated by the characters, $\Omega^+_{-n/2}(\mathbb{M}^k_{n\varpi_1})$ is identified with a quotient of a Verma module, namely, $\mathbf{M}^c(\frac{1}{4}\epsilon n,\frac{1}{2}\epsilon n)=\mathbb{M}^c(\frac{1}{4}\varepsilon n,\frac{1}{2}\varepsilon n)/(G_{-1/2}^+|\frac{1}{4}\varepsilon n,\frac{1}{2}\varepsilon n\rangle)$, called a topological Verma module\cite{Ad1,FSST}. The other one $\Omega^+_{-n/2}(\mathbb{M}^k_{-(n+2)\varpi_1})$ is obtained from $\mathbf{M}^c(\tfrac{-1}{4}\varepsilon (n+2),\tfrac{-1}{2}\varepsilon(n+2))$ by spectral flow twist $S_{-n-1}$ defined in general by 
\begin{align*}
S_\theta\colon G^{\pm}_a \mapsto G^\pm_{a\pm \theta},\quad J_a\mapsto J_a+\tfrac{c}{3}\theta\delta_{a,0},\quad L_a \mapsto L_a+\theta J_a+\tfrac{c}{6}\theta^2\delta_{n,0},\quad (\theta\in \Z),
\end{align*}
which transforms the module characters in the following way:
$$\mathrm{ch}\ S_pM(q,z)=q^{\frac{c}{6}p^2}z^{\frac{c}{3}p}\mathrm{ch}\ M(q,zq^p).$$
Therefore, the counterpart of \eqref{BGG type resolution} is 
$$0\rightarrow S_{-n-1}\mathbf{M}^c(\tfrac{-1}{4}\varepsilon (n+2),\tfrac{-1}{2}\varepsilon(n+2))\rightarrow \mathbf{M}^c(\tfrac{1}{4}\epsilon n,\tfrac{1}{2}\epsilon n)\rightarrow L_c(\tfrac{1}{4}\varepsilon n,\tfrac{1}{2}\varepsilon n)\rightarrow 0.$$

To obtain resolutions of simple $V^c(\ns_2)$-modules by Verma modules $\mathbb{M}^c(h,m)$ (they are called \emph{massive} Verma modules), we have to replace the affine Verma modules $\mathbb{M}^k_\lambda$ by thicker $V^k(\sll_2)$-modules. They are given by \emph{relaxed highest weight modules} $R_{a,b}^k$ ($a,b\in \C$) induced by the weight $\sll_2$-modules
\begin{align*}
R_{a,b}:=\mathcal{U}(\sll_2)\otimes_{\C[\Omega,h]}\C_{a,b}\simeq \bigoplus_{n>0}\C e^n v_{a,b}\oplus v_{a,b} \oplus \bigoplus_{n>0}\C f^n v_{a,b}.
\end{align*}
Here $\Omega=\frac{1}{2}h^2+ef+fe$ is the Casimir element of $\sll_2$ and $\C_{a,b}$ is the one dimensional module over $\C[\Omega,h]$ with $\Omega=a$, $h=b$. Then $R_{a,b}^k$ has the character $\mathrm{tr}_\bullet q^{L_0}z^{h/2}$
$$\mathrm{ch} R_{a,b}^k=q^{\frac{a}{2(k+2)}}z^{\frac{b}{2}}\frac{\sum_{n\in\Z}z^n}{\poch{q,zq,z^{-1}q}}=q^{\frac{a}{2(k+2)}}z^{\frac{b}{2}}\frac{\sum_{n\in\Z}z^n}{\poch{q}^3},$$
which gives the desired character on $V^c(\ns_2)$-side:
$$\mathrm{ch}\ \Omega_{-n/2}^+\left(R_{\frac{n(n+2)}{2},n}^k\right)=q^{\frac{1}{4}\varepsilon n}z^{\frac{1}{2}\varepsilon n}
\frac{\poch{-zq^{\frac{1}{2}},-z^{-1}q^{\frac{1}{2}} }}{\poch{q}^2}.$$
Indeed, we have an isomorphism $ \Omega_{-n/2}^+(R_{n(n+2)/2,n}^k)\simeq \mathbb{M}^c(\frac{1}{4}\varepsilon n,\frac{1}{2}\varepsilon n)$.
The relaxed highest weight modules and their spectral flow twists defined through
$$S_\theta\colon e_a\mapsto e_{a+\theta},\quad h_a\mapsto h_a+k \theta \delta_{a,0},\quad f_a\mapsto f_{a-\theta},\quad (\theta\in \Z),$$
are the most natural class of $V^k(\sll_2)$-modules since all the simple weight $V^k(\sll_2)$-modules are obtained\cite{Fut} as their simple quotients.
A nice realization of relaxed highest weight modules is implemented by the inverse of the quantum Drinfeld--Sokolov reduction of $V^k(\sll_2)$ introduced by Adamovi\'{c}\cite{Ad3}:
\begin{align}\label{invHR}
\begin{split}
&V^k(\sll_2)\longrightarrow \W^k(\sll_2)\otimes \Pi\\
&e(z)\mapsto 1\otimes \mathrm{e}^{\mathbf{a}}(z),\quad \tfrac{1}{2}h(z)\mapsto 1\otimes \mathbf{b}^+(z),\\
&f(z)\mapsto  (k+2)L(z)\otimes \mathrm{e}^{-\mathbf{a}}(z),\\
&\hspace{2cm}-1\otimes :\left(\mathbf{b}^{-}(z)^2+(k+1)(\partial_z\mathbf{b}^-(z) \right)\mathrm{e}^{-\mathbf{a}}(z):,
\end{split}
\end{align}
where $\Pi$ is the half-lattice vertex algebra
$$\Pi:=V_{(1+\ssqrt{-1})\Z}\underset{\pi^{(1+\ssqrt{-1})\Z}}{\otimes}\pi^{\Z\oplus \ssqrt{-1}\Z}\subset V_{\Z\oplus \ssqrt{-1}\Z},$$
and 
\begin{align*}
\mathbf{a}(z)=b_{1+\ssqrt{-1}}(z),\quad \mathbf{b}^\pm(z)=\pm\tfrac{k}{4}b_{1+\ssqrt{-1}}(z)+\tfrac{1}{2}b_{1-\ssqrt{-1}}(z).
\end{align*}
Let us introduce the following modules\cite{Feh} over $V^k(\sll_2)$ 
\begin{align*}
\mathrm{M}^k_{r,s}[\theta,\lambda]:=\mathrm{M}^{c(k)}_{r,s}\otimes \Pi_{\theta}[\lambda],\quad (r,s,\theta\in \Z,\ \lambda\in\C).
\end{align*}
Here $\mathrm{M}^{c(k)}_{r,s}$ is the Verma module of highest weight $L_0=h_{r,s}$ over the Virasoro algebra of central charge $c(k)$ with
$$c(k)=1-6\tfrac{(k+1)^2}{(k+2)},\quad h_{r,s}=\tfrac{\left (r(k+2)-s\right)^2-(k+1)^2}{4(k+2)},$$ 
and $\Pi_\theta[\lambda]$ is the $\Pi$-module given by the sum of Fock modules 
$$\Pi_\theta[\lambda]= \bigoplus_{n\in \Z}\pi^{\Z\oplus \ssqrt{-1}\Z}_{\theta \mathbf{b}^++(\lambda+n)\mathbf{a}}.$$
In particular, $\mathrm{M}_{r,s}^k[-1,\lambda]$ has the following character
\begin{align*}
\mathrm{ch}\ \mathrm{M}_{r,s}^k[-1,\lambda]=q^{h_{r,s}+\frac{k}{4}}z^{\lambda-\frac{k}{2}}\frac{\sum_{n\in\Z}z^n}{\poch{q}^3}
\end{align*}
and is isomorphic to the relaxed highest weight module
\begin{align}\label{realization of relaxed module}
M_{r,s}^k[-1,\lambda] \simeq R_{a,b}^k,\quad (a=2(k+2)(h_{r,s}+\tfrac{k}{4}),\quad b=2(\lambda-\tfrac{k}{2}))
\end{align}
under the assumption $(b-p)^2 \neq 2a+1$ for all positive odd integers $p\in\Z$. The assumption says that there is no highest weight vector in $R_{a,b}^k$ otherwise these two modules are non-isomorphic indecomposable modules. The other modules $M_{r,s}^k[\theta,\lambda]$ are obtained from $M_{r,s}^k[-1,\lambda]$ by spectral flow twists:
$$M_{r,s}^k[\theta,\lambda]\simeq S_{\theta+1}M_{r,s}^k[-1,\lambda].$$

Let $\lambda$ be generic so that \eqref{realization of relaxed module} holds. The simple quotient $\mathrm{L}_{r,s}^{c(k)}$ ($r,s\geq1$) of $\mathrm{M}_{r,s}^{c(k)}$ has a two step resolution\cite{FeiFu, IK} by Verma modules due to Feigin--Fuchs:
$$0\rightarrow \mathrm{M}^{c(k)}_{-r,s}\rightarrow \mathrm{M}^{c(k)}_{r,s}\rightarrow\mathrm{L}_{r,s}^{c(k)}\rightarrow 0.$$
Except for finitely many values of $\lambda$, $\mathrm{L}_{r,s}^{c(k)}\otimes \Pi_{-1}[\lambda]$ coincides with the simple quotient of $\mathrm{M}^k_{r,s}[-1,\lambda]$, which we denote by $\mathrm{L}^k_{r,s}[-1,\lambda]$.
In this case, the above resolution induces a resolution of $\mathrm{L}^k_{r,s}[-1,\lambda]$ by relaxed highest Verma modules
\begin{align*}
0\rightarrow \mathrm{M}_{-r,s}^k[-1,\lambda]\rightarrow \mathrm{M}_{r,s}^k[-1,\lambda]\rightarrow \mathrm{L}_{r,s}^k[-1,\lambda]\rightarrow 0.
\end{align*}
Applying the functor $\Omega^+_{-b/2}$, we obtain the following resolution of a simple $V^c(\ns_2)$-module by massive Verma modules:
\begin{align*}
0\rightarrow \mathbb{M}^c(\alpha_-,\beta_\lambda)\rightarrow \mathbb{M}^c(\alpha_+,\beta_\lambda)\rightarrow L_c(\alpha_+,\beta_\lambda)\rightarrow 0
\end{align*}
with
$$\alpha_\pm=h_{\pm r,s}+\tfrac{1}{2\varepsilon}\left(\beta_0-\beta_\lambda^2 \right),\quad \beta_\lambda=\varepsilon(\lambda+1)-1.$$
The above approach should also work at admissible levels. 
Some resolutions of $V^k(\sll_2)$-modules and $V^c(\ns_2)$-modules in this case has been obtained in Refs. \citen{KR,FSST,KS} .
It might be an interesting problem to compare these two.

We record the following dictionary\cite{FSST} on the correspondence of some classes of basic modules up to spectral flow twists. One of the aims of this article is to present a generalization for other $\W$-superalgebras, see Tab. \ref{correspondence of modules in general} at the end of this article.
\renewcommand{\arraystretch}{1.3}
\begin{table}[ph]
\caption{}
\label{correspondence of modules}
\begin{tabular}{cc}
\hline
$V^k(\sll_2)$-side & $V^c(\ns_2)$-side \\ \hline
affine Verma modules & topological Verma modules  \\
relaxed highest weight modules & massive Verma modules \\ \hline
\end{tabular}
\end{table}
\renewcommand{\arraystretch}{1}
\subsection{Fusion rules: rational case}

Let us compare the fusion rings for $L_k(\sll_2)$-mod and $L_c(\ns_2)$-mod when $k$ is a non-negative integer, i.e. the \emph{rational} case, when both categories are semisimple braided tensor categories with finitely many simple objects.
Although the equivalence Theorem \ref{equivalence of categories} is stated at the level of abelian categories, it is strong enough to deduce the relation of their fusion rules in this case by the theory of simple current extensions in the language of braided tensor category\cite{CKL,CKM,YY,CGNS}. 

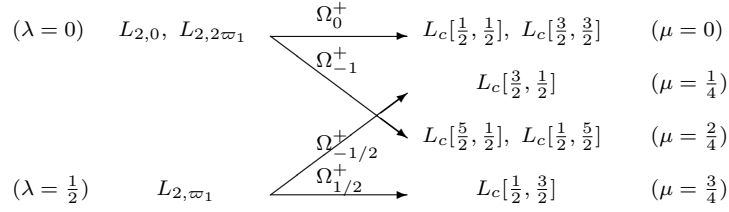
\begin{figure}[htbp]
\centering
\caption{Correspondence of simple modules for $k=2$}
\setlength{\unitlength}{1mm}
\begin{picture}(83,30)(0,7)\label{block-wise equiv}
\put(-4,28){\footnotesize($\lambda=0$)}
\put(-4,7){\footnotesize($\lambda=\tfrac{1}{2}$)}
\put(10,28){\footnotesize$L_{2,0},\ L_{2,2\varpi_1}$}
\put(15,7){\footnotesize$L_{2,\varpi_1}$}
\put(50,28){\footnotesize $L_c[\tfrac{1}{2},\tfrac{1}{2}],\ L_c[\tfrac{3}{2},\tfrac{3}{2}$]}
\put(57,21){\footnotesize$L_c[\tfrac{3}{2},\tfrac{1}{2}]$}
\put(50,14){\footnotesize$L_c[\tfrac{5}{2},\tfrac{1}{2}],\ L_c[\tfrac{1}{2},\tfrac{5}{2}]$}
\put(57,7){\footnotesize$L_c[\tfrac{1}{2},\tfrac{3}{2}]$}
\put(80,28){\footnotesize($\mu=0$)}
\put(80,21){\footnotesize($\mu=\tfrac{1}{4}$)}
\put(80,14){\footnotesize($\mu=\tfrac{2}{4}$)}
\put(80,7){\footnotesize($\mu=\tfrac{3}{4}$)}
\put(30,28){\vector(1,0){18}}
\put(30,28){\vector(4,-3){18}}
\put(30,7){\vector(4,3){18}}
\put(30,7){\vector(3,0){18}}
\put(36,30){\footnotesize$\Omega^+_{0}$}
\put(36,24){\footnotesize$\Omega^+_{-1}$}
\put(36,13){\footnotesize$\Omega^+_{-1/2}$}
\put(36,8.5){\footnotesize$\Omega^+_{1/2}$}
\end{picture}
\end{figure}

Recall that the list of simple $L_k(\sll_2)$-modules is given by 
$$L_{k,0},\ L_{k,\varpi_1},\ \cdots,\ L_{k,k\varpi_1}$$
and that the list\cite{Ad1} of simple $L_c(\ns_2)$-modules is given by the highest weight modules
\begin{align*}
L_{c}[u,v]:=L_c(h_{u,v},s_{u,v}),\quad h_{u,v}=\tfrac{1}{k+2}(uv-\tfrac{1}{4}),\ s_{u,v}=\tfrac{u-v}{k+2},
\end{align*}
where $\mathbf{u}=(u,v)$ runs through the set
$$u,v\in \tfrac{1}{2}+\Z_{\geq0},\ 0\leq u,v,u+v\leq k+2.$$
They consist of all the unitary minimal representations\cite{DVPYZ} of $\ns_{2,c}$.
The module categories decompose into
\begin{align*}
L_k(\sll_2)\text{-mod}=\bigoplus_{i\in\Z_2}L_k(\sll_2)\text{-mod}^{[\frac{i}{2}]},\quad L_c(\ns_2)\text{-mod}=\bigoplus_{j\in\Z_{k+2}}L_c(\ns_2)\text{-mod}^{[\frac{j}{k+2}]}
\end{align*}
and the lists of simple modules lying in each block are given by
\begin{align*}
L_{k,p\varpi_1}\ (p\equiv i\ \text{mod}\ 2),\quad L_c[\mathbf{u}]\ (u-v\equiv j\  \text{mod}\ k+2).
\end{align*}
In Fig. \ref{block-wise equiv}, we illustrate the block-wise equivalence for $k=2$.

The fusion rules are given by
$$L_{k,p\varpi_1}\boxtimes L_{k,q\varpi_1}\simeq \bigoplus_{r} L_{k,r\varpi_1},\quad L_c[\mathbf{u}_1]\boxtimes  L_c[\mathbf{u}_2]\simeq \bigoplus_{\mathbf{u}_3}L_c[\mathbf{u}_3] $$
where $r$ runs through
\begin{itemize}
\item $|p-q|\leq r\leq \mathrm{Min}\{p+q,2k-(p+q)\},$\\
$r\equiv p+q \text{ mod } 2,$
\end{itemize}
and $\mathbf{u}_3$ through\cite{Ad2} either of 
\begin{itemize}
\item $|\mathbf{u}_1^+-\mathbf{u}_2^+|<\mathbf{u}_3^+<\mathrm{Min}\{\mathbf{u}_1^++\mathbf{u}_2^+,2(k+2)-(\mathbf{u}_1^++\mathbf{u}_2^+)\}$,\\
$\mathbf{u}_3^-=\mathbf{u}_1^-+\mathbf{u}_2^-$,
\item $|\mathbf{u}_1^+-\mathbf{u}_2^+|<(k+2)-\mathbf{u}_3^+<\mathrm{Min}\{\mathbf{u}_1^++\mathbf{u}_2^+,2(k+2)-(\mathbf{u}_1^++\mathbf{u}_2^+)\}$,\\
$\mathbf{u}_3^-=\mathbf{u}_1^-+\mathbf{u}_2^-\pm (k+2)$,
\end{itemize}
where we set $\mathbf{u}^\pm_i=u_i\pm v_i$.
In particular, $L_{k,k\varpi_1}$ is a simple current and satisfies
$$L_{k,k\varpi_1}\boxtimes  L_{k,p\varpi_1}\simeq L_{k,(k-p)\varpi_1}.$$

For a more transparent comparison of the fusion rings, we use a general phenomenon of simple current extensions by lattice theories\cite{YY,CKM,CGNS} applied for the decomposition 
\begin{align*}
L_c(\ns_2)\otimes V_{\ssqrt{-1}\Z}\simeq L_{k,0}\otimes V_{\ssqrt{-2(k+2)}\Z}\oplus L_{k,\varpi_1}\otimes V_{\frac{(k+2)}{\ssqrt{-2(k+2)}}+\ssqrt{-2(k+2)}\Z}.
\end{align*}
Let $V$ be a rational $C_2$-cofinite vertex operator algebra so that the module category $V$-mod is a semisimple braided tensor category of finitely many simple objects. We assume that $V$-mod has simple current modules $S_{a}$ which form an abelian group $G$ by fusion product
$$S_a\boxtimes S_b\simeq S_{a+b}.$$
 
Note that the fusion ring $\mathcal{K}(V)$ of $V$-modules is naturally a ring over $\Z[G]$ by fusion product $a.M=S_a \boxtimes M$ and a ring graded by the dual group $G^\vee=\mathrm{Hom}(G,\C^\times)$ by monodromy: $\mathcal{M}_{S_a,M}=\xi(a)\mathrm{id}_{S_a\boxtimes M}$ with $\xi\in G^\vee$. By a simple current extension of $V$ by a lattice theory, we mean a vertex operator superalgebra extension of the form
$$\mathcal{E}=\bigoplus_{a\in G}S_a\otimes V_{F(a)+L}$$
for some non-degenerate integral lattice $L$ with a group homomorphism $F \colon G\hookrightarrow L'/L$. 
Here $L'$ is the dual lattice $L'=\{a\in \Q\otimes_\Z L\mid (a,L)\subset \Z\}$ with which we may parametrize the simple $V_L$-modules $V_{\lambda+L}$ ($\lambda\in L'/L$). Take the sublattice $N\subset L'$ so that $F\colon G\simeq N/L\subset L'/L$.
Then the fusion rings $\mathcal{K}(V)$ and $\mathcal{K}(\mathcal{E})$ are related by the formulas
\begin{align*}
\mathcal{K}(\mathcal{E})\simeq \left(\mathcal{K}(V)\underset{\Z[N/L]}{\otimes} \Z[L'/L] \right)^{N/L},\quad \mathcal{K}(V)\simeq \left(\mathcal{K}(\mathcal{E})\underset{\Z[N'/L]}{\otimes} \Z[L'/L] \right)^{N'/L}.
\end{align*}
In our setting, they are
\begin{align}
&\label{SCA}\mathcal{K}(L_c(\ns_2))\simeq \left(\mathcal{K}(L_k(\sll_2))\underset{\Z[\Z_2]}{\otimes}\Z[\Z_{2(k+2)}]\right)^{\Z_2},\\
&\mathcal{K}(L_k(\sll_2))\simeq \left(\mathcal{K}(L_c(\ns_2))\underset{\Z[\Z_{k+2}]}{\otimes}\Z[\Z_{2(k+2)}]\right)^{\Z_{k+2}}.
\end{align}
We note that in Ref. \cite{Ad2} the fusion rules of $L_c(\ns_2)$ were actually derived from \eqref{SCA}.

Finally, we compare the fusion rings from the common ground, that is the fusion ring of their Heisenberg coset, known as the parafermion algebra.
For $k\in \Z_{\geq0}$, it enjoys the \emph{level-rank duality}\cite{ALY} and is identified with a simple principal $\W$-algebra 
\begin{align}\label{relation of KS}
\Com\left(\pi^h,L_k(\sll_2)\right)\simeq \W_{\ell}(\sll_{k}),\quad \ell=-k+\frac{k+2}{k+1}.
\end{align} 
The $\W$-algebra $\W_{\ell}(\sll_{k})$ is rational and $C_2$-cofinite\cite{Ar1,Ar2} and known as one of the discrete series representation. The module category $\W_\ell(\sll_k)$-mod is a semisimple braided tensor category with simple modules parametrized by those of $L_2(\sll_k)$, say $\mathbf{L}_\W^k(\lambda)$, satisfying the same fusion rules.\cite{AvE1,FKW,Cr}
On the other hand, the double commutant $\Com(\W_\ell(\sll_k),L_k(\sll_2))$ extending the Heisenberg vertex algebra $\pi^h$ turns out to be a lattice vertex algebra $V_{\sqrt{2k}\Z}$. Then $L_k(\sll_2)$ decomposes as a module over $\W_\ell(\sll_k)\otimes V_{\sqrt{2k}\Z}$ into 
\begin{align*}
L_k(\sll_2)\simeq \bigoplus_{i\in \Z_k}\mathbf{L}_{\W}^k(2\varpi_i)\otimes V_{\frac{2i}{\sqrt{2k}}+\sqrt{2k}\Z}
\end{align*}
and, similarly, 
\begin{align*}
L_c(\ns_2)\simeq \bigoplus_{i \in \Z_k}\mathbf{L}_{\W}^k(2 \varpi_i)\otimes V_{\frac{(k+2)i}{\sqrt{(k+2)k}}+\sqrt{(k+2)k}\Z }.
\end{align*}
Therefore, $L_k(\sll_2)$ and $L_c(\ns_2)$ are simple current extensions by lattice theories, which imply the following description of fusion rings:
\begin{align*}
&\mathcal{K}(L_k(\sll_2))\simeq \left(\mathcal{K}(L_2(\sll_k))\underset{\Z[\Z_k]}{\otimes}\Z[\Z_{2k}] \right)^{\Z_k},\\
&\mathcal{K}(L_c(\sll_2))\simeq \left(\mathcal{K}(L_2(\sll_k))\underset{\Z[\Z_{(k+2)k}]}{\otimes}\Z[\Z_{(k+2)k}] \right)^{\Z_k}.
\end{align*}

\section{Feigin--Semikhatov conjecture}
\subsection{Main statements}
To obtain a fruitful generalization of the relation between $V^k(\sll_2)$ and $V^c(\ns_2)$, it is better to view $V^c(\ns_2)$ as the principal $\W$-superalgebra associated with $\sll_{2|1}$:
\begin{align*}
V^c(\ns_2)\simeq \W^\ell(\sll_{2|1}),\quad c=-3(2\ell+1),
\end{align*} 
and then rewrite the relation \eqref{relation of KS} into
\begin{align}
(k+2)(\ell+1)=1.
\end{align}
Note that the numbers two and one in each factor are the dual Coxeter numbers of $\sll_2$ and $\sll_{2|1}$. This implies that  the relation between $V^k(\sll_2)$ and $V^c(\ns_2)$ is a variant of the celebrated \emph{Feigin--Frenkel duality}\cite{FF} for the principal $\W$-algebras, which asserts an isomorphism
\begin{align}\label{Feigin Frenkel duality}
\W^k(\sll_n)\simeq \W^{\check{k}}(\sll_n),\quad (k+n)(\check{k}+n)=1.
\end{align}
in the case of type A.

Recall that the $\W$-superalgebras in general are defined from the affine vertex superalgebras $V^\kappa(\g)$ for (simple) basic classical Lie superalgebras $\g$ with even supersymmetric invariant bilinear forms $\kappa$ via the quantum Drinfeld--Sokolov reductions parametrized by the conjugacy classes of even nilpotent elements $f$, 
\begin{align*}
\W^\kappa(\g,f):=H^0_{\mathrm{DS},f}(V^\kappa(\g)).
\end{align*}
When $\g=\sll_{n|m}$ ($n\neq m$), the bilinear forms $\kappa$ are all proportional to the super trace $\mathrm{str}$ and we identify $\kappa$ with its ratio $k\in \C$. The conjugacy classes are in one-to-one correspondence with the pair of partitions of $n$ and $m$. 
For $\g=\sll_2$ they are $2,1^2$ and correspond to the Virasoro vertex algebra and $V^k(\sll_2)$, respectively. For $\g=\sll_{2|1}$, they are $(2|1), (1^2|1)$ and correspond to the $\mathcal{N}=2$ superconformal algebra and $V^k(\sll_{2|1})$, respectively.

The set of all the even nilpotent elements forms an algebraic subvariety $\mathcal{N}$ inside $\sll_{n|m}$ and decomposes into conjugacy classes $\mathcal{N}=\sqcup_\lambda \mathcal{N}_\lambda$. 
There is a unique class, called regular or \emph{principal}, characterized as maximizing the dimension $\mathrm{dim}\ \mathcal{N}_\lambda$.
The partition is given by $(n|m)$ and the corresponding $\W$-superalgebra $\W^k(\sll_{n|m})=\W^k(\sll_{n|m},f_{\mathrm{prin}})$ with $f_{\mathrm{prin}}=f_{n|m}$ is called principal. 
If $\g=\sll_n$, then there is a unique class, called \emph{subregular}, of the second largest dimension. The partition is given by $(n-1)^11^1$ and the corresponding $\W$-algebra $\W^k(\sll_n,f_{\mathrm{sub}})$ with $f_{\mathrm{sub}}=f_{(n-1)^11^1}$ is called subregular.

Feigin and Semikhatov\cite{FS} conjectured that the pair $V^k(\sll_2)$ and $\W^\ell(\sll_{2|1})$ generalizes to the series of pairs
\begin{align}\label{Feigin Semikhatov pair}
\W^k(\sll_{n},f_{\mathrm{sub}}),\quad \W^\ell(\sll_{n|1},f_{\mathrm{prin}}).
\end{align}

The $\W$-superalgebras \eqref{Feigin Semikhatov pair} have rank one Heisenberg vertex subalgebras as the maximal affine subalgebras: we normalize the generators 
\begin{align}
H_+(z)=\sum_{n\in\Z}H_{+,n}z^{-n-1},\quad H_-(z)=\sum_{n\in\Z}H_{-,n}z^{-n-1}
\end{align}
so that the eigenvalues of the zero mode $H_{\pm,0}$ on the $\W$-superalgebras are normalized to $\Z$. Then the OPEs are
\begin{align*}
H_+(z)H_+(w)\sim \frac{\varepsilon_+}{(z-w)^2},\quad H_-(z)H_-(w)\sim \frac{\varepsilon_-}{(z-w)^2}.
\end{align*}
with 
\begin{align}\label{heisenberg}
\varepsilon_+=\tfrac{n-1}{n}(k+n)-1,\quad \varepsilon_-=-\tfrac{n}{n-1}(\ell+n-1)+1
\end{align}
The precise analogy of the Feigin--Frenkel duality for the pairs \eqref{Feigin Semikhatov pair} is the duality 
between the Heisenberg cosets of these $\W$-superalgebras:
\begin{theorem}\cite{CL1,CGN}\label{Feigin-Semikhatov}
There is an isomorphism of vertex algebras 
\[{\bf FS}\colon\Com \left( \pi^{H_+}, \W^{k} (\mathfrak{sl}_n, f_{\mathrm{sub}}) \right) \simeq \Com \left( \pi^{H_-}, \W^{\ell} (\mathfrak{sl}_{n|1},f_{\mathrm{prin}}) \right)\]
for $(k,\ell)\neq(-n+\frac{n}{n-1},-(n-1)+\frac{n-1}{n})$ satisfying the relation 
\begin{align}\label{duality level}
(k+n)(\ell+n-1)=1.
\end{align}
\end{theorem}
These pairs indeed enjoy the coset constructions:
\begin{theorem}\cite{CGN}\label{Kazama-Suzuki}
There are isomorphisms of vertex superalgebras
\begin{align*}
&{\bf KS}\colon\W^{\ell}(\mathfrak{sl}_{n|1},f_{\mathrm{prin}})\xrightarrow{\simeq}\Com(\pi^{H_+^{\Delta}},\W^k(\mathfrak{sl}_n,f_{\mathrm{sub}})\otimes V_{\Z}),\\[2mm]
&{\bf FST}\colon\W^k(\mathfrak{sl}_n,f_{\mathrm{sub}})\xrightarrow{\simeq}\Com(\pi^{H_-^{\Delta}},\W^{\ell}(\mathfrak{sl}_{n|1},f_{\mathrm{prin}})\otimes V_{\ssqrt{-1}\Z})
\end{align*}
for all $(k,\ell)$ satisfying \eqref{duality level}. Here we set 
$$H_+^{\Delta}(z)=-H_+(z)+b_1(z),\quad H_-^{\Delta}(z)=H_-(z)+b_{\ssqrt{-1}}(z).$$
\end{theorem}
Note that we can prove Theorem \ref{Kazama-Suzuki} for $n=2$  for $V^k(\sll_2)$ by checking the OPEs directly. However, we cannot use this approach literally since we do not know the whole defining OPEs of $\W^k(\sll_n,f_{\mathrm{sub}})$ and $\W^\ell(\sll_{n|1},f_{\mathrm{prin}})$. 
There are two approaches to overcome this difficulty: one is the \emph{uniqueness} property of these algebras which is am important feature for the hook-type $\W$-superalgebras in general as proved by Creutzig and Linshaw \cite{CL1,CL2} and the other one is \emph{free field realizations} of these $\W$-superalgebras. 

\subsection{Uniqueness property} 
The $\W$-algebra $\W^k(\sll_n,f_{\mathrm{sub}})$ has a free strong generating set of the form
$$\W^k(\sll_n,f_{\mathrm{sub}})=\W\left(1,2,\cdots,n-1,(\tfrac{n}{2})^2\right),$$
that is, a generating set whose conformal weights are $1,2\cdots, \tfrac{n}{2}$ with described multiplicity such that the ordered monomials of their negative modes give a basis of the algebra. The field of conformal weight one is the Heisenberg field $H_+(z)$ and the field of conformal weight two is the Virasoro field. The fields of conformal weight $\frac{n}{2}$ are those generators which have nontrivial $H_{+,0}$-weights, namely $\pm1$.

Similarly, the $\W$-superalgebra $\W^\ell(\sll_{n|1},f_{\mathrm{rpin}})$ has a free strong generating set of the form
$$\W^\ell(\sll_{n|1},f_{\mathrm{prin}})=\W\left(1,2,\cdots,n,(\tfrac{n+1}{2})^2\right).$$
Again, the field of conformal weight one is the Heisenberg field $H_-(z)$, the field of conformal weight two is the Virasoro field and the fields of conformal weight $\frac{n+1}{2}$ are those generators which have nontrivial $H_{-,0}$-weights, namely $\pm1$ and are the only odd generators.

The Heisenberg cosets 
\begin{align}\label{Heisnenberg coset}
\Com \left( \pi^{H_+}, \W^{k} (\mathfrak{sl}_n, f_{\mathrm{sub}}) \right),\quad \Com \left( \pi^{H_-}, \W^{\ell} (\mathfrak{sl}_{n|1},f_{\mathrm{prin}}) \right)
\end{align}
have a special feature: they have a common strong generating set of the form
\begin{align}\label{generating type}
\W(2,3,\cdots,2n+1)
\end{align}
which is the same as the principal $\W$-algebra $\W^k(\sll_{2n+1})$, in which case it is free. 
In our case, the generating set is not free and indeed has two linearly independent relation at conformal dimension $2n+4$, which is a consequence of a $q$-character \eqref{q-character of coset} below.
 
The important fact is that one-parameter families of vertex algebras of type $\W(2,3,\cdots,N)$ are always obtained as quotients of the two-parameter \emph{universal $\W_\infty$-algebra} $\W_\infty[c,\lambda]$ constructed mathematically by Linshaw \cite{L}, whose existence has been conjectured for decades, see the references therein. It is a vertex algebra over the commutative ring $\C[c,\lambda]$ and has a strong free generating set of the form
$$\W_\infty[c,\lambda]=\W(2,3,4,\cdots).$$
By setting $W_m(z)$ to be the generating field of conformal degree $m$, $c$ parametrizes the central charge of the Virasoro field $W_2(z)$ and $\lambda$ parametrizes the undetermined scalar appearing in the OPE
\begin{align*}
&W_2(z)W_5(w)\sim \frac{-5(16(c+2)\lambda-37) W_3(w)}{(z-w)^4}\\
&\hspace{4cm}+\frac{-(16(c+2)\lambda-55)\partial_w W_3(w)}{(z-w)^3}+\frac{5 W_5(w)}{(z-w)^2}+\frac{\partial_w W_5(w)}{(z-w)}.
\end{align*}
The remaining OPEs are uniquely determined by the Jacobi identities under the assumption that there is an automorphism sending $W_{m}(z)\mapsto (-1)^m W_{m}(z)$.
Both of the coset algebras \eqref{Heisnenberg coset} and $\W^k(\sll_n)$ ($n\geq3$) are obtained as one-parameter families of vertex algebras from $\W_\infty[c,\lambda]$ with different algebraic relations between $c$ and $\lambda$, called \emph{truncation curves} in general.
In the case $\W^k(\sll_n)$, the parameters $c, \lambda$ are 
\begin{align*}
&c=(n-1)\left(1-\frac{n(n+1)(k+n-1)^2}{k+n}\right),\\
&\lambda=-\frac{n+k}{(n-2)(n^2+nk-n-2)(n^2+nk +n+2k)},
\end{align*}
giving a coordinate of the truncation curve $\mathbf{C}_{\mathrm{pr}}(n)$
$$\lambda(n-2)(3n^2-n-2+c(n+2))=(n^2-1).$$
The $\W$-algebra $\W^k(\sll_n)$ is obtained from the specialization of $\W_\infty[c,\lambda]$ at these values via the quotient by the ideal generated by the null  vector of conformal weight $n+1$. 
On the other hand, \eqref{Heisnenberg coset} is obtained from $\W_\infty[c,\lambda]$ by the specialization 
\begin{align*}
&c=-\frac{n (n(n-3)+k(n-2)+1) (n(n-2)+k(n-1)-1)}{k + n},\\
&\lambda=-\frac{(k+n)(k+n-1)}{(n(n-4)+k (n-3)+2) (n(n-2)+k(n-1)-2) (n^2+nk+k)},
\end{align*}
which parametrizes the truncation curve $\mathbf{C}_{\mathrm{sub}}(n)$
\begin{align*}
\begin{split}
&c^3 \lambda^2 (n-3) (n-2) (n+1) (n+2)\\
&+c^2 \lambda (-3 (4 \lambda+1)+\lambda n (8 n^4-33 n^3+2 n^2+99 n-28)-2 (n-1) n (n^2-n-5))\\
&+c n \binom{13 \lambda+\lambda^2 (n-1) (12 n^4-88 n^3+179 n^2+\lambda n-20)\hspace{3cm}}{\hspace{3cm}-\lambda n (10 n^3-37 n^2+18 n+34)+n (n^2-2n-1)+2}\\
&-n^2 (\lambda (n-5) (3 n+1)-2 (n+1)) (5 \lambda+(4 \lambda+1) (n-3) n+2)\\
&=0.
\end{split}
\end{align*}

The above uniqueness ensures the Heisenberg cosets \eqref{Heisnenberg coset} coincide. 
The whole $\W$-superalgebra $\W^k(\sll_n,f_{\mathrm{sub}})$ (resp.\ $\W^\ell(\sll_{n|1},f_{\mathrm{prin}})$) is characterized\cite{CL1} by the vertex superalgebra containing the tensor product of the coset and $\pi^{H_+}$ (resp.\ $\pi^{H_-}$), extended by even (resp.\ odd) highest weight vectors of conformal weight $\frac{n}{2}$ (resp.\ $\frac{n+1}{2}$)  and of Heisenberg weight $\pm1$.
This gives a proof of Theorem \ref{Kazama-Suzuki} at generic levels.

The level-rank duality for the parafermion algebra \eqref{relation of KS} can be systematically found\cite{CL1} at the intersection points of the truncation curves $\mathbf{C}_{\mathrm{sub}}(n)$ and $\mathbf{C}_{\mathrm{pr}}(m)$.
\begin{theorem}\cite{ALY,ACL2,CL1}\label{level-rank duality} There is an isomorphism of vertex algebras 
$$\Com\left(\pi^{H_+},\W_k(\sll_n,f_{\mathrm{sub}})\right)\simeq \W_r(\sll_m)$$
for the following shifted levels $(k+n,r+m)$ with $m\geq3$:
\begin{align*}
&\left(\frac{n+m}{n-1},\frac{m+1}{m+n}\right),\quad \left(\frac{n}{n+m-1},\frac{m-1}{m+n-1}\right),\quad \left(\frac{n-m}{n-m-1},\frac{n-m}{n-m-1}\right).
\end{align*}
\end{theorem}
The condition $m\geq3$ comes from the construction of $\W_\infty[c,\lambda]$ that $W_3(z)$ does not vanish, which prohibits degenerations to the Virasoro vertex algebra.
The case $m=2$ is subtle: the statement is true in the first case \cite{CGNS}, but not in the second case \cite{ACGY}. The latter case is of particular interest since the Heisenberg coset is the singlet algebra, which is an simple current extension of the Virasoro vertex algebra of infinite order. The algebra $\W_k(\sll_n,f_{\mathrm{sub}})$ itself is identified with the chiral algebras of Argyres--Douglas theories of type $(A_1,A_{2n-1})$\cite{ACGY}.

\subsection{A digression on $q$-characters}
Let us briefly discuss the $q$-character of the Heisenberg cosets \eqref{Heisnenberg coset}. 
Since the character of $\W^k(\sll_n,f_{\mathrm{sub}})$ is expressed as 
\begin{align*}
\mathrm{tr}_{\W^k(\sll_n,f_{\mathrm{sub}})}(q^{L_0}z^{H_{+,0}})=\frac{1}{(q,\cdots q^{n-1},zq^{\frac{n}{2}}z^{-1}q^{\frac{n}{2}};q)_\infty},
\end{align*} 
the $q$-character of the Heisenberg cosets is given by the residue
\begin{align}\label{character via residue}
F(q):=\int \frac{dz}{z} \frac{1}{(q^2,\cdots q^{n-1},zq^{\frac{n}{2}},z^{-1}q^{\frac{n}{2}};q)_\infty}.
\end{align}
As a $z$-independent form, it can be rewritten into
\begin{align}\label{q-character of coset}
\frac{1}{(q;q)_\infty^2(q,\cdots,q^n;q)_\infty}\sum_{m\in\Z}(-1)^m\Phi_{m}(q)
\left(q^{\frac{n+1}{2}+m},q^{\frac{n+1}{2}-m};q\right)_\infty q^{\frac{m(m+1)}{2}}
\end{align}
by the unary false theta functions
\begin{align}\label{some function}
\Phi_m(q)=\sum_{N=0}^\infty(-1)^Nq^{\frac{N(N+1)}{2}+Nm}.
\end{align}
From this, we find that there are two independent linear relations among the PBW basis for strong generators at conformal dimension $2n+4$ in \eqref{generating type}. 

There is a combinatorial meaning\cite{GR}, called \emph{melting crystals}, for \eqref{character via residue} when multiplied by the character of the rank one Heisenberg vertex algebra $\pi$, which is $1/(q;q)_\infty $.  
It is based on the following fact.
When we extend $\W_{\infty}[c,\lambda]$ by $\pi$, the resulting vertex algebra $\W_{1+\infty}=\pi\otimes \W_\infty[c,\lambda]$ is of type $\W(1,2,3,\cdots)$ and the $q$-character is given by the MacMahon function
$$\frac{1}{\prod_{n=1}^\infty(1-q^n)^n}.$$
It is known as the generating function counting plane partitions (i.e., three dimensional partitions located in the region $\{x,y,z\geq 0\}\subset \mathbb{R}^3$). In the same line, the function $F(q)/(q;q)_\infty$ is the generating function of plane partitions with a pit\cite{BFM} at $(2,n+1)$ as in Fig. \ref{plane partition with a pit}.
The first excluded plane partition is the one in Fig. \ref{excluded partition}, which has $2(n+1)$ boxes. This corresponds to the fact that the field $W_{2n+2}(z)$ drops off from the strong generating set in \eqref{generating type}.

\setlength{\unitlength}{1mm}
\begin{figure}[htbp]
\begin{minipage}[b]{6.3cm}
\centering
\caption{Plane partiton with a pit at $(2,n+1)$}
\setlength{\unitlength}{1mm}
\begin{picture}(30,26)(0,0)\label{plane partition with a pit}
\put(0,0){\line(1,0){24}}
\put(0,4){\line(1,0){20}}
\put(0,8){\line(1,0){20}}
\put(0,12){\line(1,0){20}}
\put(0,16){\line(1,0){20}}
\put(0,20){\line(1,0){8}}
\put(0,0){\line(0,1){24}}
\put(4,0){\line(0,1){24}}
\put(8,0){\line(0,1){24}}
\put(12,0){\line(0,1){16}}
\put(16,0){\line(0,1){16}}
\put(4,20){\line(1,-1){8}}
\put(4,16){\line(1,-1){4}}
\put(4,16){\line(1,1){4}}
\put(4,12){\line(1,1){4}}
\put(8,12){\line(1,1){4}}
\put(12,16){\line(1,-1){4}}
\put(12,12){\line(1,1){4}}
\put(1.5,1){\footnotesize7}
\put(5.5,1){\footnotesize5}
\put(9.5,1){\footnotesize2}
\put(1.5,9){\footnotesize4}
\put(5.5,9){\footnotesize2}
\put(1.5,13){\footnotesize3}
\put(1.5,17){\footnotesize1}
\put(9.5,9){\footnotesize1}
\put(13.5,9){\footnotesize0}
\put(13.5,1){\footnotesize1}
\put(1.5,4.5){\footnotesize$\vdots$}
\put(1.5,20.5){\footnotesize$\vdots$}
\put(5.5,4.5){\footnotesize$\vdots$}
\put(5.5,20.5){\footnotesize$\vdots$}
\put(9.5,4.5){\footnotesize$\vdots$}
\put(13.5,4.5){\footnotesize$\vdots$}
\put(17.5,1){\footnotesize$\cdots$}
\put(17.5,5){\footnotesize$\cdots$}
\put(17.5,9){\footnotesize$\cdots$}
\put(17.5,13){\footnotesize$\cdots$}
\put(20,7){$\left. \begin{array}{c}\\ \\ \\ \end{array} \right\}$}
\put(28,7){\footnotesize $n+1$}
\end{picture}
\end{minipage}
\begin{minipage}[b]{6.2cm}
\centering
\caption{First excluded partition}
\setlength{\unitlength}{1mm}
\begin{picture}(15,26)(0,0)\label{excluded partition}
\put(0,0){\line(1,0){20}}
\put(0,0){\line(0,1){24}}
\put(0,4){\line(1,0){8}}
\put(0,8){\line(1,0){8}}
\put(0,12){\line(1,0){8}}
\put(0,16){\line(1,0){8}}
\put(4,0){\line(0,1){16}}
\put(8,0){\line(0,1){16}}
\put(1.5,1){\footnotesize1}
\put(5.5,1){\footnotesize1}
\put(1.5,9){\footnotesize1}
\put(5.5,9){\footnotesize1}
\put(1.5,13){\footnotesize1}
\put(5.5,13){\footnotesize1}
\put(1.5,4.5){\footnotesize$\vdots$}
\put(5.5,4.5){\footnotesize$\vdots$}
\put(9,7){$\left. \begin{array}{c}\\ \\ \\ \end{array} \right\}$}
\put(17,7){\footnotesize $n+1$}
\end{picture}
\end{minipage}
\end{figure}

\subsection{Free field realization}

Beside the uniqueness property of the $\W$-superalgebras based on $\W_{\infty}[c,\lambda]$, the free field realization of  $\W$-superalgebras is also a powerful tool as in the original proof\cite{FF} of the Feigin--Frenkel duality \eqref{Feigin Frenkel duality}. 
In general, we have the Miura map \cite{KRW}
$$\W^k(\g,f)\hookrightarrow V^{k'}(\g_0)\otimes SB(\g_{1/2}),$$
which embeds the $\W$-superalgebra into a tensor product of the affine vertex algebra $V^{k'}(\g_0)$ at a suitable level and the symplectic boson associated with the symplectic vector superspace $\g_{1/2}$. Here we use the good grading $\g=\bigoplus_{j\in \frac{1}{2}\Z}\g_{j}$ used in the quantum Drinfeld--Sokolov reduction.
In our cases,
\begin{align}\label{Miura}
\W^k(\sll_n,f_{\mathrm{sub}})\hookrightarrow V^{k'}(\sll_2)\otimes \pi^{k+n}_{\h^\perp},\quad 
\W^\ell(\sll_{n|1},f_{\mathrm{prin}})\hookrightarrow V^{\ell'}(\gl_{1|1})\otimes \pi^{\ell+n-1}_{\h^\perp}.
\end{align}
Here $\h^\perp\subset\g$ is the subspace of the standard Cartan subalgebra $\h$ orthogonal to $\sll_2$ and $\gl_{1|1}$ respectively.  
When $k,\ell$ are generic, the images coincide with the joint kernels of screening operators associated with highest weight vectors of Fock modules for the Heisenberg vertex algebras. In the weighted Dynkin diagrams below, they are canonically attached to the nodes of weight one.

\begin{figure}[htbp]
\begin{minipage}[b]{6.5cm}
\centering
\caption{$\g=\sll_n$}
\setlength{\unitlength}{1mm}
\begin{picture}(41,13)(0,5)\label{sln}
\put(0,10){\circle{2}}
\put(-1,13){\footnotesize$1$}
\put(-1,5){\footnotesize$\alpha_1$}
\put(1,10.3){\line(1,0){8}}
\put(10,10){\circle{2}}
\put(9,13){\footnotesize$1$}
\put(9,5){\footnotesize$\alpha_2$}
\put(11,10.3){\line(1,0){6}}
\put(18.5,9.4){$\cdot$}
\put(20,9.4){$\cdot$}
\put(21.5,9.4){$\cdot$}
\put(24.5,10.3){\line(1,0){6}}
\put(32,10){\circle{2}}
\put(31,13){\footnotesize$1$}
\put(31,5){\footnotesize$\alpha_{n-2}$}
\put(33,10.3){\line(1,0){8}}
\put(42,10){\circle{2}}
\put(41,13){\footnotesize$0$}
\put(41,5){\footnotesize$\alpha_{n-1}$}
\end{picture}
\end{minipage}
\begin{minipage}[b]{5cm}
\centering
\caption{$\g=\sll_{n|1}$}
\setlength{\unitlength}{1mm}
\begin{picture}(41,13)(60,5)\label{ssln}
\put(60,10){\circle{2}}
\put(59,13){\footnotesize$1$}
\put(59,5){\footnotesize$\alpha_1$}
\put(61,10.3){\line(1,0){8}}
\put(70,10){\circle{2}}
\put(69,13){\footnotesize$1$}
\put(69,5){\footnotesize$\alpha_2$}
\put(71,10.3){\line(1,0){6}}
\put(78.5,9.4){$\cdot$}
\put(80,9.4){$\cdot$}
\put(81.5,9.4){$\cdot$}
\put(84.5,10.3){\line(1,0){6}}
\put(92,10){\circle{2}}
\put(91,13){\footnotesize$1$}
\put(91,5){\footnotesize$\alpha_{n-1}$}
\put(93,10.3){\line(1,0){8}}
\put(102,10){\circle{2}}
\put(100.7,9.3){\footnotesize$\times$}
\put(101,13){\footnotesize$0$}
\put(101,5){\footnotesize$\alpha_{n}$}
\end{picture}
\end{minipage}
\end{figure}

We compose them with the Wakimoto realization
\begin{align}\label{Wakimoto}
V^{k'}(\sll_2)\hookrightarrow \beta\gamma\otimes \pi^{k+n}_\h,\quad V^{\ell'}(\gl_{1|1})\hookrightarrow bc\otimes \pi_\h^{\ell+n-1} \simeq V_{\Z}\otimes \pi_\h^{\ell+n-1}
\end{align}
and then the bosonization\cite{FMS} of the $\beta\gamma$-system 
\begin{align}\label{FMS}
\beta\gamma\hookrightarrow \Pi:=V_{(1+\ssqrt{-1})\Z}\underset{\pi^{(1+\ssqrt{-1})\Z}}{\otimes}\pi^{\Z\oplus \ssqrt{-1}\Z}.
\end{align}
In the end, we arrive at the following free field realizations of the $\W$-superalgebras:
\begin{align}\label{FFR}
& \W^k(\sll_n,f_{\mathrm{sub}})\hookrightarrow \Pi\otimes \pi^{k+n}_{\h},\quad \W^\ell(\sll_{n|1},f_{\mathrm{prin}})\hookrightarrow V_\Z\otimes \pi^{\ell+n-1}_{\h}.
\end{align}
At each step \eqref{Wakimoto}-\eqref{FMS}, the images are described by the kernel of a single screeing operator. Therefore, at generic levels $k,\ell$, the images of \eqref{FFR} are described by the joint kernels of screening operators. 
The screening operators for \eqref{Wakimoto} are attached to the nodes of weight zero in Fig. \ref{sln}-\ref{ssln}, respectively. The screening operator for \eqref{FMS} adds one \emph{odd} node to Fig. \ref{sln} as below. 
\begin{figure}[htbp]
\caption{Screening operators for $\W^k(\sll_n,f_{\mathrm{sub}})$}
\centering
\setlength{\unitlength}{1mm}
\begin{picture}(53,13)(-1,5)\label{scr for subreg}
\put(0,10){\circle{2}}
\put(-1,13){\footnotesize$1$}
\put(-1,5){\footnotesize$\alpha_1$}
\put(1,10.3){\line(1,0){8}}
\put(10,10){\circle{2}}
\put(9,13){\footnotesize$1$}
\put(9,5){\footnotesize$\alpha_2$}
\put(11,10.3){\line(1,0){6}}
\put(18.5,9.4){$\cdot$}
\put(20,9.4){$\cdot$}
\put(21.5,9.4){$\cdot$}
\put(24.5,10.3){\line(1,0){6}}
\put(32,10){\circle{2}}
\put(31,13){\footnotesize$1$}
\put(31,5){\footnotesize$\alpha_{n-2}$}
\put(33,10.3){\line(1,0){8}}
\put(42,10){\circle{2}}
\put(41,13){\footnotesize$0$}
\put(41,5){\footnotesize$\alpha_{n-1}$}
\put(43,10.3){\line(1,0){8}}
\put(52,10){\circle{2}}
\put(50.7,9.3){\footnotesize$\times$}
\end{picture}
\end{figure}

Now, we may apply the Feigin--Frenkel duality for the Virasoro vertex algebras for each (even) node in Fig. \ref{ssln}-\ref{scr for subreg} since the Heisenberg cosets concentrate on the Heisenberg vertex algebras by \eqref{FFR}:
\begin{align}\label{FFR for cosets}
\begin{split}
&\Com\left(\pi^{H_+},\W^k(\sll_{n},f_{\mathrm{sub}})\right)\hookrightarrow\Com\left(\pi^{H_+},\pi^{\Z\oplus \ssqrt{-1}\Z}\otimes \pi^{k+n}_\h\right),\\
&\Com\left(\pi^{H_-},\W^\ell(\sll_{n|1},f_{\mathrm{prin}})\right)\hookrightarrow\Com\left(\pi^{H_-},\pi^{\Z}\otimes \pi^{\ell+n-1}_\h\right).
\end{split}
\end{align}
This proves Theorem \ref{Feigin-Semikhatov} for generic levels.
Note that these cosets form continuous families of subspaces inside the common Heisenberg vertex algebra and that the Feigin--Frenkel duality for the Virasoro vertex algebras asserts that they are the same as vector spaces. Then it follows that the isomorphism in Theorem \ref{Feigin-Semikhatov} remains true for all levels except for the level $(k,\ell)=(-n+\frac{n}{n-1},-(n-1)+\frac{n-1}{n})$, when $\pi^{H_\pm}$ degenerate to commutative vertex algebras.
Theorem \ref{Kazama-Suzuki} can be obtained in the same way by using \eqref{FFR} itself.

\begin{remark}
The screening realizations of the Heisenberg cosets \eqref{FFR for cosets} coincide with those proposed in Refs. \citen{PR,GR2} based on Miura operators. 
Therefore, they coincide with the vertex algebras introduced by Bershtein--Feigin--Merzon \cite{BFM}, which act on the cohomology of the moduli spaces of spiked instantons of Nekrasov.\cite{RSYY}
\end{remark}

From the perspective of the free field realization \eqref{FFR}, the inverse Drinfeld--Sokolov reduction \eqref{invHR} is naturally explained: the appearance of the Virasoro vertex algebra comes from its free free field realization inside $\pi^{k+2}_\h$, whose image is also characterized by an even screening operator. 
We stress that the screening operators for $\W^k(\sll_2)$ and $V^k(\sll_2)$ have different form and that an automorphism of $\pi^{k+2}_\h\otimes \Pi$ is used to identify them.
Nontheless, it is generalized to our setting together with the $\W^\ell(\sll_{n|1},f_{\mathrm{prin}})$-side in the following way:
\begin{align}\label{invHR in general}
\xymatrix@=18pt{\W^k(\sll_n,f_{\mathrm{sub}}) \ar[rr] \ar@{}[d]|{\bigcap} & & \W^\ell(\sll_{n|1},f_{\mathrm{prin}})\otimes V_{\ssqrt{-1}\Z} \ar@{}[d]|{\bigcap} \\
\W^k(\sll_n)\otimes \Pi \ar[rr] \ar@{}[d]|{\bigcap}& &  \W^{\ell-1}(\sll_n) \otimes \pi^{h_\perp}\otimes V_{\Z\oplus \ssqrt{-1}\Z} \ar@{}[d]|{\bigcap}\\
\pi^{k+n}_\h\otimes \Pi  \ar[rr]&& \pi^{\ell+n-1}_\h\otimes V_{\Z\oplus \ssqrt{-1}\Z}.
 }
\end{align}
Here $h_\perp(z)$ is the Heisenberg field corresponding to the matrix $\sum_{i=1}^nE_{i,i}+nE_{n+1,n+1}$ inside  $\sll_{n|1}$ in terms of the elementary matrices $E_{i,j}$. 
The embedding 
$$\W^k(\sll_n,f_{\mathrm{sub}})\hookrightarrow \W^k(\sll_n)\otimes \Pi$$
has been introduced and used by Adamovi\'{c}--Kawasetsu--Ridout\cite{AKR} for the Bershadsky--Polyakov algebra $\W^k(\sll_3,f_{\mathrm{sub}})$ and by Fehily\cite{Feh} in general, to construct relaxed highest weight modules which are a source of logarithmic conformal field theories. It is obvious that we may convert the horizontal arrow in \eqref{invHR in general} by using the other coset construction (Theorem \ref{Kazama-Suzuki}) by adding necessary lattice vertex superalgebras. In the end, we find that these two $\W$-superalgebras are contained as 
\begin{align*}
\W^k(\sll_n,f_{\mathrm{sub}})\subset \W^k(\sll_n)\otimes V_{\Z\oplus \ssqrt{-1}\Z}\supset \W^\ell(\sll_{n|1},f_{\mathrm{prin}}),
\end{align*} 
where they just stretch out in different lattices, namely $(1+\ssqrt{-1})\Z$ and $\Z$. 
This explains why the module categories of $\W^k(\sll_n,f_{\mathrm{sub}})$ and $\W^\ell(\sll_{n|1},f_{\mathrm{prin}})$ look similar in general.

\section{Correspondence of representation categories}

Since we have the coset construction (Theorem \ref{Kazama-Suzuki}) as in the case of $V^k(\sll_2)$ and $V^c(\ns_2)$ in Section 2.1, we can generalize the equivalence of weight module categories and compare the fusion data in the rational setting.

\subsection{Equivalence of categories}

By using the normalized Heisenberg fields in \eqref{heisenberg}, we introduce the category of weight modules consisting of all the modules which decompose into direct sums of Fock modules when restricted to the Heisenberg vertex subalgebras, that is \emph{Kazhdan--Lusztig objects} for $\widehat{\gl}_1$.
Motivated by this, we write 
$$\mathbf{KL}_k^{n-1,1},\quad \mathbf{KL}_\ell^{n|1}$$
for the weight module categories for $\W^k(\sll_n,f_{\mathrm{sub}})$ and $\W^\ell(\sll_{n|1},f_{\mathrm{prin}})$, respectively. 
They decompose into blocks 
\begin{align*}
\mathbf{KL}_k^{n-1,1}=\bigoplus_{[\lambda]\in \C/\Z}\mathbf{KL}_{k,[\lambda]}^{n-1,1},\quad
\mathbf{KL}_\ell^{n|1}=\bigoplus_{[\mu]\in \C/\Z} \mathbf{KL}_{\ell,[\mu]}^{n|1}
\end{align*}
in terms of the eigenvalues of $H_{\pm,0}$ as the spectrum of $H_{\pm,0}$ of the $\W$-superalgebras is $\Z$.
Then we have functors
\begin{align*}
\Omega^+_\xi\colon \mathbf{KL}_{k,[\lambda]}^{n-1,1}\rightarrow \mathbf{KL}_{\ell,[(\varepsilon_--1)\xi ]}^{n|1},\quad
\Omega^-_{\xi'}\colon \mathbf{KL}_{\ell,[\mu]}^{n|1}\rightarrow \mathbf{KL}_{k,[(\varepsilon_++1)\xi' ]}^{n-1,1},
\end{align*}
which assign the multiplicity space $\Omega^+_\xi(M)$ of the modules $\pi^{H_+^\Delta}_\xi$ inside $M\otimes V_\Z$ to a $\W^k(\sll_n,f_{\mathrm{sub}})$-module $M$ and the multiplicity space $\Omega^-_{\xi'}(N)$ of the modules $\pi^{H_-^\Delta}_{\xi'}$ inside $N\otimes V_{\ssqrt{-1}\Z}$ to a $\W^\ell(\sll_{n|1},f_{\mathrm{prin}})$-module $N$, respectively. 
Note that the ratios $\varepsilon_\pm\pm1$ are the norms of $H_\pm^{\Delta}(z)$.
\begin{theorem}\cite{CGNS}\label{equiv of category in general}
The functors 
\begin{align}
\Omega_{-\lambda}^+\colon \mathbf{KL}_{k,[\lambda]}^{n-1,1}\ \rightleftarrows  \ \mathbf{KL}_{\ell,[(1-\varepsilon_-)\lambda]}^{n|1} \colon \Omega^-_{(1-\varepsilon_-)\lambda}
\end{align}
are quasi-inverse to each other and thus give an equivalence of categories. 
\end{theorem}
Thanks to the relation \eqref{duality level}, the value $1-\varepsilon_-$ coincides with the ratio of the Heisenberg fields $H_{\pm}(z)$, namely $\varepsilon_-/\varepsilon_+$. Later, we will use the square root:
$$\epsilon=\sqrt{1-\varepsilon_-}=\sqrt{\varepsilon_-/\varepsilon_+}.$$

\subsection{Fusion rules: rational case}

Here we compare the fusion data of the full subcategories 
$$\W_k(\sll_n,f_{\mathrm{sub}})\text{-mod}\subset \mathbf{KL}_k^{n-1,1},\quad \W_\ell(\sll_{n|1},f_{\mathrm{prin}})\text{-mod}\subset \mathbf{KL}_\ell^{n|1} $$
when $(k,f_{\mathrm{sub}})$ form an exceptional pair\cite{AvE2,KW}, i.e.
\begin{align}\label{exceptional level}
k=-n+\tfrac{n+r}{n-1},\quad (r\geq2,\ \mathrm{gcd}(n+r,n-1)=1).
\end{align}
In this case, $\W_k(\sll_n,f_{\mathrm{sub}})$ is rational \cite{AvE2} and $C_2$-cofinite \cite{Ar2}, and thus $\W_k(\sll_n,f_{\mathrm{sub}})$-mod is a semisimple braided tensor category with finitely many simple objects. These algebraic and categorical properties are inherited to $\W_\ell(\sll_{n|1},f_{\mathrm{prin}})$-side \cite{CGN,CGNS}.

To describe the fusion data for the simple exceptional $\W$-algebra $\W_k(\sll_n,f_{\mathrm{sub}})$, it is better to use the level-rank duality (Theorem \ref{level-rank duality}):
$$\Com\left(\pi^{H_+},\W_{-n+\frac{n+r}{n-1}}(\sll_n,f_{\mathrm{sub}})\right) \simeq \W_{-r+\frac{r+n}{r+1}}(\sll_r).$$
The module category $\W_{-r+\frac{r+n}{r+1}}(\sll_r)$-mod is a semisimple braided tensor category with simple modules $\mathbf{L}_\W^r(\lambda)$ satisfying the same fusion rules as $L_n(\sll_r)$-modules:
\begin{align}
\mathcal{K}\left(\W_{-r+\frac{r+n}{r+1}}(\sll_r)\right)\xrightarrow{\simeq} \mathcal{K}(L_n(\sll_r)),\quad \mathbf{L}_\W^r(\lambda)\mapsto L_{n,\lambda}.
\end{align}
Note that $\mathbf{L}_\W^r(n\varpi_i)$ ($0\leq i\leq r-1$) are simple currents. They form a group $\Z_r$ since 
$$\mathbf{L}_\W^r(n\varpi_i) \boxtimes \mathbf{L}_\W^r(n\varpi_j)\simeq \mathbf{L}_\W^r(n\varpi_{i+j}).$$
Then the fusion ring is a ring over $\Z[\Z_r]$ by fusion products
$$\mathbf{L}_{\W}^r(n\varpi_i)\boxtimes \mathbf{L}_\W^r(\lambda)\simeq \mathbf{L}_\W^r(\sigma^i\lambda)$$
where $\sigma$ permutes the fundamental weights $\sigma \varpi_i=\varpi_{i+1}$,
and also admits a grading by the dual group $\Z_r^\vee=\mathrm{Hom}(\Z_r,\C^\times)$ coming from the monodromy
$$\mathcal{M}_{\mathbf{L}_\W^r(n \varpi_i),\mathbf{L}_\W^r(\lambda)}=\zeta_r^{-i \pi_{P/Q}(\lambda)}$$
with $\zeta_r=\mathrm{exp}(2\pi\ssqrt{-1}/r)$ and $\pi_{P/Q}(\varpi_i)=i$.
Then $\W_k(\sll_n,f_{\mathrm{sub}})$ decomposes\cite{CL1} into 
\begin{align*}
\W_k(\sll_n,f_{\mathrm{sub}})\simeq \bigoplus_{i\in\Z_r} \mathbf{L}_\W^r(n\varpi_i)\otimes V_{\frac{ni}{\sqrt{nr}}+\sqrt{nr}\Z}.
\end{align*}
Therefore, $\W_k(\sll_n,f_{\mathrm{sub}})$ is a simple current extension of $\W_{-r+\frac{r+n}{r+1}}(\sll_r)$ by a lattice theory, which implies that the fusion ring is described\cite{AvE2,CGNS} as 
\begin{align*}
\mathcal{K}(\W_k(\sll_n,f_{\mathrm{sub}}))
\simeq 
\left(\mathcal{K}\left(\W_{-r+\frac{r+n}{r+1}}(\sll_r)\right) \underset{\Z[\Z_{r}]}{\otimes}\Z[\Z_{nr}]\right)^{\Z_{r}}
\simeq \mathcal{K}(L_r(\sll_n))
\end{align*}
by the level-rank duality between the fusion rings of affine vertex algebras\cite{Fr,OS,CGNS}.

Accordingly, we have\cite{CGNS}
\begin{align*}
\W_\ell(\sll_{n|1},f_{\mathrm{prin}})\simeq \bigoplus_{i\in \Z_r}\mathbf{L}_\W^r(n \varpi_i)\otimes V_{\frac{(n+r)i}{\sqrt{(n+r)r}}+\sqrt{(n+r)r}\Z}
\end{align*}
and then 
\begin{align*}
\mathcal{K}(\W_\ell(\sll_{n|1},f_{\mathrm{prin}}))\simeq \left(\mathcal{K}(L_n(\sll_r))\underset{\Z[\Z_r]}{\otimes}\Z[\Z_{(n+r)r}] \right)^{\Z_r}.
\end{align*}
The fusion data is  compared in the following way:
\begin{theorem}\cite{CGNS}
For $k=-n+\frac{n+r}{n-1}$ and $\ell=-(n-1)+\frac{n-1}{n+r}$ as in \eqref{exceptional level}, we have isomorphisms of rings
\begin{align*}
&\mathcal{K}(\W_\ell(\sll_{n|1},f_{\mathrm{prin}}))\simeq \left(\mathcal{K}(\W_k(\sll_n,f_{\mathrm{sub}}))\underset{\Z[\Z_n]}{\otimes}\Z[\Z_{n(n+r)}]\right)^{\Z_n},\\
&\mathcal{K}(\W_k(\sll_n,f_{\mathrm{sub}}))\simeq \left(\mathcal{K}(\W_\ell(\sll_{n|1},f_{\mathrm{prin}}))\underset{\Z[\Z_{n+r}]}{\otimes}\Z[\Z_{n(n+r)}]\right)^{\Z_{n+r}},
\end{align*}
and $\mathcal{K}(\W_k(\sll_n,f_{\mathrm{sub}}))\simeq \mathcal{K}(L_r(\sll_n))$.
\end{theorem}

\section{Correspondence of intertwining operators}
\subsection{Gluing approach}
Let us discuss the monoidal correspondence beyond the rational case. 
In this generality, we do not know whether there is a suitable tensor structure on the category of weight modules and, even if it exists, the technique of tensor category at hand seems not to be sufficient. Thus, let us compare the bare spaces of logarithmic intertwining operators. 
For this purpose, it is meaningful to take a closer look at the coset construction (Theorem \ref{Feigin-Semikhatov}). 
Decompose the $\W$-superalgebras $\W^k(\sll_n,f_{\mathrm{sub}})$ and $\W^\ell(\sll_{n|1},f_{\mathrm{prin}})$ as modules over the tensor product of the Heisenberg coset \eqref{Heisnenberg coset} and the Heisenberg vertex algebras $\pi^{H_\pm}$ inside them:
\begin{align}
\W^k(\sll_n,f_{\mathrm{sub}})\simeq \bigoplus_{a\in \Z}\mathscr{C}^k_{+,a}\otimes \pi^{H_+}_a,\quad 
\W^\ell(\sll_{n|1},f_{\mathrm{prin}})\simeq \bigoplus_{a\in \Z}\mathscr{C}^\ell_{-,a}\otimes \pi^{H_-}_a.
\end{align}
Here $\mathscr{C}^k_{+,a}$ and $\mathscr{C}^\ell_{-,a}$ are the space of highest weight vector of the Fock modules.
Then Theorem \ref{Feigin-Semikhatov} implies that 
$$\mathscr{C}^k_{+,a}\simeq \mathscr{C}^\ell_{-,a}$$
as modules over the Heisenberg cosets. In other words, the difference between these two $\W$-superalgebras comes from the norms of the generators of the Heisenberg vertex subalgebras $H_+(z)$ and $H_-(z)$. This difference can be seen most explicitly when their simple quotients $\W_k(\sll_n,f_{\mathrm{sub}})$ and $\W_\ell(\sll_{n|1},f_{\mathrm{prin}})$ degenerate to free field algebras:
\begin{align*}
\W_{-2m+\frac{2m+1}{2m-1}}(\sll_{2m},f_{\mathrm{sub}})\simeq V_{\sqrt{2m}\Z},\quad \W_{-(2m-1)+\frac{2m-1}{2m+1}}(\sll_{2m|1},f_{\mathrm{prin}})\simeq V_{\sqrt{2m+1}\Z}.
\end{align*}

The difference coming from the normalizations of Heisenberg fields is indeed responsible for the difference of module categories: the block-wise equivalence of categories is a formulation to get rid of this effect and the difference of the periodicity of simple currents (equivalently, spectral flows) in the rational case manifests this difference. Now, it is natural to seek for a way to change the normalizations of Heisenberg fields more directly.
Note that it is equivalent to swap the Fock modules $\pi^{H_+}_a$ and $\pi^{H_-}_a$ simultaneously.

The very role is played by the \emph{relative semi-infinite cohomology functor} founded by Feigin\cite{Fei} and Frenkel--Garland--Zuckerman \cite{FGZ} as conjectured by Creutzig--Linsahw \cite{CL1,CL2}.
In our case, we use the the semi-infinite cohomology of $\widehat{\mathfrak{gl}}_1$ relative to the horizontal subalgebra $\gl_1$. 
Let $\pi^{A^+}$ and $\pi^{A^-}$ be Heisenberg vertex algebras generated by the fields $A^+(z)$ and $A^-(z)$ satisfying the OPEs
$$A^+(z)A(w)^+\sim \frac{1}{(z-w)^2},\quad A^-(z)A^-(w)\sim \frac{-1}{(z-w)^2}.$$
Then the diagonal one 
\begin{align}
\pi^{\mathrm{diag}}:=\langle A^+(z)\otimes 1+1\otimes A^-(z)\rangle \subset\pi^{A^+}\otimes \pi^{A^-}
\end{align}
is a commutative vertex algebra which acts on the tensor product of Fock modules $\pi_\alpha^{A^+}\otimes \pi^{A^-}_\beta$ of $\pi^{A^+}\otimes \pi^{A^-}$. In this setting, we have 
\begin{align}\label{relcoh}
H^{\mathrm{rel},n}_\infty \left(\gl_1; \pi^{A^+}_\alpha\otimes \pi^{A^-}_{\beta}\right)\simeq \delta_{\alpha+\beta,0}\delta_{n,0}\C[\mathrm{e}^\alpha \otimes \mathrm{e}^\beta].
\end{align}

Therefore, to replace the Fock module $\pi^{H_+}_a$ inside $\W^k(\sll_n,f_{\mathrm{sub}})$ with $\pi^{H_-}_a$ simultaneously in order to obtain the shape of $\W^\ell(\sll_{n|1},f_{\mathrm{prin}})$, we introduce the following module
\begin{align}\label{gluing1}
K_{+\rightarrow -}:=\bigoplus_{a\in \Z}\pi_{-a}^{\ssqrt{-1}H_+}\otimes \pi^{H_-}_a
\end{align}
over $\pi^{\ssqrt{-1}H_+}\otimes \pi^{H_-}$. Then \eqref{relcoh} already implies that 
\begin{align}\label{kernel gluing}
H_\infty^{\mathrm{rel},0}\left(\gl_1; \W^k(\sll_n,f_{\mathrm{sub}})\otimes K_{+\rightarrow-}\right)\simeq \bigoplus_{a\in \Z}\mathscr{C}_{+,a}^k\otimes \pi_a^{H_-}\simeq \W^\ell(\sll_{n|1},f_{\mathrm{prin}})
\end{align}
as modules over $\Com(\pi^{H_-}, \W^\ell(\sll_{n|1},f_{\mathrm{prin}}))\otimes \pi^{H_-}$. 
It is crucial here to observe that \eqref{gluing1} admits a vertex superalgebra structure 
$$K_{+\rightarrow -}\simeq V_\Z\otimes \pi^{A^-}$$
by changing basis of the two dimensional subspace of Heisenberg fields inside. 
Here again appears the lattice vertex superalgebra $V_\Z$ used in the coset construction (Theorem \ref{Kazama-Suzuki}). Thanks to this, on can show that \eqref{kernel gluing} is an isomorphism of vertex superalgebras. The appearance of $V_\Z$ is not a mere coincidence: indeed, when we manipulate the opposite direction in the same way, we find that $V_{\ssqrt{-1}\Z}$ appears as the main part of the gluing object
\begin{align*}
&\W^k(\sll_n,f_{\mathrm{sub}})\simeq H_\infty^{\mathrm{rel},0}\left(\gl_1;\W^\ell(\sll_{n|1},f_{\mathrm{prin}})\otimes K_{-\rightarrow +}\right),\\ 
&K_{-\rightarrow+}\simeq V_{\ssqrt{-1}\Z}\otimes \pi^{A^+}.
\end{align*}
To summarize, we obtain the following reformulation of the coset construction:
\begin{theorem}\cite{CGNS}\label{reformulation by gluing}
There exist isomorphisms of vertex superalgebras
\begin{align}
\begin{split}
&\W^\ell(\sll_{n|1},f_{\mathrm{prin}})\simeq H_\infty^{\mathrm{rel},0}\left(\gl_1; \W^k(\sll_n,f_{\mathrm{sub}})\otimes K_{+\rightarrow-}\right),\\
&\W^k(\sll_n,f_{\mathrm{sub}})\simeq H_\infty^{\mathrm{rel},0}\left(\gl_1;\W^\ell(\sll_{n|1},f_{\mathrm{prin}})\otimes K_{-\rightarrow +}\right),
\end{split}
\end{align}
with 
\begin{align}\label{gluing objects}
K_{+\rightarrow -}=V_\Z\otimes \pi^{A^-},\quad K_{-\rightarrow+}=V_{\ssqrt{-1}\Z}\otimes\pi^{A^+}.
\end{align}
\end{theorem}

\subsection{Correspondence}

To study the correspondence of intertwining operators, let us start with reformulating the equivalence of weight module categories (Theorem \ref{equiv of category in general}) through the gluing approach. For the coset functor, we have used the multiplicity space of Fock modules with non-trivial highest weights. There is a similar room for the gluing objects \eqref{gluing objects}: we may replace the rank one Heisenberg vertex algebra by the Fock modules. We set
\begin{align}\label{gluing modules}
\begin{split}
&K_{+\rightarrow-}^\lambda:=V_\Z\otimes \pi^{A^-}_\lambda\simeq \bigoplus_{\mu\in\Z}\pi^{\ssqrt{-1}H_+}_{-\mu-\epsilon^{-1}\lambda}\otimes \pi^{H_-}_{\mu+\epsilon \lambda},\\
&K_{-\rightarrow+}^{\lambda}:=V_{\ssqrt{-1}\Z}\otimes \pi^{A^+}_{\lambda}\simeq \bigoplus_{\mu\in \Z}\pi^{\ssqrt{-1}H_-}_{-\mu-\epsilon \lambda}\otimes \pi^{H_+}_{\mu+\epsilon^{-1}\lambda}.
\end{split}
\end{align} 
By Theorem \ref{reformulation by gluing}, they give rise to functors
\begin{align*}
&H_{+,\lambda}^{\mathrm{rel}}(\bullet):=H_{\infty}^{\mathrm{rel},0}(\gl_1; \bullet\otimes K_{+\rightarrow-}^\lambda)\colon \mathbf{KL}_k^{n-1,1}\rightarrow  \mathbf{KL}_\ell^{n|1},\\
&H_{-,\lambda}^{\mathrm{rel}}(\bullet):=H_{\infty}^{\mathrm{rel},0}(\gl_1; \bullet\otimes K_{-\rightarrow+}^\lambda)\colon \mathbf{KL}_\ell^{n|1}\rightarrow \mathbf{KL}_k^{n-1,1}.
\end{align*}
It turns out that they are equivalent to the coset functors:
\begin{theorem}\cite{CGNS}
We have natural isomorphisms 
\begin{align*}
&H_{+,\epsilon\lambda}^{\mathrm{rel}}\simeq \Omega_{-\lambda}^+\colon \mathbf{KL}_{k,[\lambda]}^{n-1,1}\rightarrow \mathbf{KL}_{\ell,[\epsilon^2\lambda]}^{n|1},\\
&H_{-,\epsilon\lambda}^{\mathrm{rel}}\simeq \Omega_{\epsilon^2\lambda}^-\colon \mathbf{KL}_{\ell,[\epsilon^2\lambda]}^{n|1} \rightarrow  \mathbf{KL}_{k,[\lambda]}^{n-1,1}.
\end{align*}
In particular, the functors $H_{+,\epsilon\lambda}^{\mathrm{rel}}$ and $H_{-,\epsilon\lambda}^{\mathrm{rel}}$ give an equivalence of categories 
\begin{align}\label{equivalence of categories via relcoh} 
H_{+,\epsilon\lambda}^{\mathrm{rel}}\colon \mathbf{KL}_{k,[\lambda]}^{n-1,1}\rightleftarrows \mathbf{KL}_{\ell,[\epsilon^2\lambda]}^{n|1} \colon H_{-,\epsilon\lambda}^{\mathrm{rel}}.
\end{align}
\end{theorem}
Now, it is clear how to incorporate with \eqref{equivalence of categories via relcoh} the spaces of intertwining operators. 
Returning back to \eqref{gluing modules}, we just use the intertwining operators of Fock modules.

More concretely, we fix the base of the space of intertwining operators as follows:
\begin{align*}
I_{\pi^{A^+}}\binom{\pi^{A^+}_{\lambda_3}}{\pi^{A^+}_{\lambda_1}\ \pi^{A^+}_{\lambda_2}}=\C Y^+(\cdot,z),\quad I_{\pi^{A^-}}\binom{\pi^{A^-}_{\lambda_3}}{\pi^{A^-}_{\lambda_1}\ \pi^{A^-}_{\lambda_2}}=\C Y^-(\cdot,z)
\end{align*}
for $\lambda_3=\lambda_1+\lambda_2$ so that the values of highest weight vectors are given by
\begin{align*}
Y^\pm(\mathrm{e}^{\lambda_1},z)\mathrm{e}^{\lambda_2}=z^{\pm\lambda_1\lambda_2}\mathrm{exp}
\left(\lambda_1\sum_{m\geq0}\frac{A^\pm_{-m}}{m}z^m \right)\mathrm{e}^{\lambda_3}.
\end{align*}
Then, given a logarithmic intertwining operator of $\W^k(\sll_n,f_{\mathrm{sub}})$-modules
\begin{align*}
\mathcal{Y}(\bullet,z)\colon M_1\otimes M_2\rightarrow M_3\{z\}[\log z],
\end{align*}
we have a logarithmic intertwining operator 
\begin{align*}
\mathcal{Y}(\bullet,z)\otimes Y^-(\bullet,z)\colon (M_1\otimes \pi^{A^-}_{\epsilon \lambda_1})\otimes (M_2\otimes \pi^{A^-}_{\epsilon \lambda_2})\rightarrow (M_3 \otimes \pi^{A^-}_{\epsilon \lambda_3})\{z\}[\log z],
\end{align*}
which naturally extends to the whole complex defining the relative semi-infinite cohomology. It is straightforward to show that it descends to an intertwining operator at the level of cohomology which we denote by
\begin{align*}
\mathbb{H}_+(\mathcal{Y})\colon H_{+,\epsilon\lambda_1}^{\mathrm{rel}}(M_1)\otimes H_{+,\epsilon\lambda_2}^{\mathrm{rel}}(M_2)\rightarrow H_{+,\epsilon\lambda_3}^{\mathrm{rel}}(M_3)\{z\}[\log z].
\end{align*}
In this way, we obtain a linear map $\mathbb{H}_+$ which assigns $\W^\ell(\sll_{n|1},f_{\mathrm{prin}})$-intertwining operators to  $\W^k(\sll_{n},f_{\mathrm{sub}})$-intertwining operators. Obviously, we also have a linear map in the opposite direction, which we denote by $\mathbb{H}_-$. 
\begin{theorem}\cite{CGNS}
The linear maps $\mathbb{H}_\pm$ are isomorphisms between the spaces of (logarithmic) intertwining operators
\begin{align*}
&\mathbb{H}_+\colon I\binom{M_3}{M_1\ M_2}\overset{\sim}{\rightarrow} I\binom{H_{+,\epsilon\lambda_3}^{\mathrm{rel}}(M_3)}{H_{+,\epsilon\lambda_1}^{\mathrm{rel}}(M_1)\ H_{+,\epsilon\lambda_2}^{\mathrm{rel}}(M_2)},\\
&\mathbb{H}_-\colon I\binom{N_3}{N_1\ N_2}\overset{\sim}{\rightarrow} I\binom{H_{-,\epsilon\mu_3}^{\mathrm{rel}}(N_3)}{H_{-,\epsilon\mu_1}^{\mathrm{rel}}(N_1)\ H_{+,\epsilon\mu_2}^{\mathrm{rel}}(N_2)}.
\end{align*}
\end{theorem}
\subsection{Some basic modules and resolutions}
We present some basic modules and resolutions generalizing \S \ref{Some basic modules and resolutions}. 
The results here will be contained in Ref. \cite{CLN} which treats the case of hook-type $\W$-superalgebras\cite{CL1,CL2} in Feigin--Frenkel type duality more generally.

Let the level $k$ be irrational. 
The $\W^k(\sll_n,f_{\mathrm{sub}})$-modules playing the role of the (contragradient of) affine Verma modules for $V^k(\sll_2)$ are \emph{Wakimoto modules} \cite{G}, which are defined through the composition of \eqref{Miura} and \eqref{Wakimoto}:
\begin{align}\label{Wakimoto realization of subregular} 
\W^k(\sll_n,f_{\mathrm{sub}})\hookrightarrow \beta\gamma\otimes \pi^{k+n}_\h.
\end{align}
This embedding is actually induced from the Wakimoto realization \cite{FF88,Frenkel} of the affine vertex algebra $V^k(\sll_n)$ 
\begin{align}\label{wakimoto realization}
 V^k(\sll_n)\hookrightarrow \mathcal{A}_{\mathfrak{n}_+}\otimes \pi^{k+n}_\h,\quad \mathcal{A}_{\mathfrak{n}_+}:=\beta\gamma^{\otimes \dim \mathfrak{n}_+},
\end{align}
which can be extended to a resolution 
\begin{align}\label{BGG for vacuum}
0\rightarrow V^k(\sll_n)\rightarrow D^k_{0}\rightarrow D^k_{1}\rightarrow \cdots \rightarrow 0
\end{align}
by Wakimoto modules of $V^k(\sll_n)$:
\begin{align*}
D^k_{p}=\bigoplus_{l(w)=p}\mathbb{W}^k_{w^{-1}*0},\quad \mathbb{W}^k_\lambda:=\mathcal{A}_{\mathfrak{n}^+}\otimes \pi^{k+n}_{\h,\lambda},\quad (\lambda\in \h^*).
\end{align*}
The sum runs over the elements $w$ of Weyl group of $\sll_n$ whose length are $p$ and $*$ denotes the dot action $w*\lambda=w(\lambda+\rho)-\rho$ with $\rho=\sum_{i}\varpi_i$ the Weyl vector. For Weyl modules  $\mathbb{V}_\lambda^k$ ($\lambda\in P_+$), we have similar resolutions
\begin{align}\label{BGG for weyl module}
0\rightarrow \mathbb{V}^k_\lambda\rightarrow D^k_{\lambda,0}\rightarrow D^k_{\lambda,1}\rightarrow \cdots \rightarrow 0
\end{align}
by replacing $D^k_{p}$ with $D^k_{\lambda,p}=\bigoplus_{l(w)=p}\mathbb{W}^k_{w^{-1}*\lambda}$. 
Introduce the following $\W^k(\sll_n,f_{\mathrm{sub}})$-modules:
$$T_\lambda^{k,+}=H_{\mathrm{DS},f_{\mathrm{sub}}}^0(\mathbb{V}_\lambda^k),\quad \mathbb{W}^{k,+}_{\lambda}:=\beta\gamma\otimes \pi^{k+n}_{\h,\lambda}.$$
The latter is called a Wakimoto module of $\W^k(\sll_n,f_{\mathrm{sub}})$ and $\W^k(\sll_n,f_{\mathrm{sub}})$ acts on it via \eqref{Wakimoto realization of subregular}.
Then, by applying $H_{\mathrm{DS},f_{\mathrm{sub}}}^0$ to \eqref{BGG for weyl module}, we obtain a resolution
\begin{align}\label{resolution of subregular}
0\rightarrow T_\lambda^{k,+}\rightarrow E^k_{\lambda,0}\rightarrow E^k_{\lambda,1}\rightarrow \cdots \rightarrow 0.
\end{align}
with 
$E^k_{\lambda,p}\simeq \bigoplus_{l(w)=p}\mathbb{W}^{k,+}_{w^{-1}*\lambda}$. 
When $\g=\sll_2$, this resolution for $\lambda=n\varpi_1$ is the contragradient dual of \eqref{BGG type resolution}.
Since \eqref{Wakimoto realization of subregular} maps 
\begin{align*}
&L(z)\mapsto (\tfrac{n}{2}\partial:\beta(z)\gamma(z):-:\gamma(z)\partial\beta(z):)\otimes 1\\
&\hspace{2cm}+1\otimes (L_{\mathrm{sug}}(z)+(1-\tfrac{1}{k+n})\partial \rho(z)-\tfrac{n}{2}\partial \varpi_{n-1}(z))\\
&H_+(z)\mapsto -:\beta(z)\gamma(z):\otimes1+1\otimes \varpi_{n-1}(z),
\end{align*}
with $L_{\mathrm{sug}}(z)$ the Virasoro field by the Segal--Sugarawa construction, 
the character $\mathrm{tr}_\bullet z^{H_{+,0}}q^{L_0}$ of the Wakimoto module $\mathbb{W}^{k,+}_{\lambda}$ is given by
\begin{align*}
\mathrm{ch}\ \mathbb{W}^{k,+}_{\lambda}=\frac{q^{\Delta_\lambda^k}(q^\frac{n}{2}z)^{(\lambda,\varpi_{n-1})}}{\poch{q}^{n-1}\poch{zq^{\frac{n}{2}},z^{-1}q^{1-\frac{n}{2}}}}
\end{align*}
with 
$\Delta_{\lambda}^k=\tfrac{1}{2(k+n)}(\lambda,\lambda)-(1-\tfrac{1}{k+n})(\rho,\lambda).$
Note that $\Delta_{\lambda}^k$ agrees with the conformal weight of the free field module $\pi^{k+n}_{\h,\lambda}$ of $\W^k(\sll_n)$ as implicated by \eqref{invHR in general}.

Introduce the following $\W^\ell(\sll_{n|1},f_{\mathrm{prin}})$-modules:
$$T^{\ell,-}_\lambda=H_{+,\epsilon(\lambda,\varpi_{n-1})}^{\mathrm{rel}}(T^{k,+}_\lambda),\quad W^\ell_\lambda=H_{+,\epsilon(\lambda,\varpi_{n-1})}^{\mathrm{rel}}(\mathbb{W}^{k,+}_\lambda).$$
Applying the functor $H_{+,\epsilon(\lambda,\varpi_{n-1})}^{\mathrm{rel}}$ to \eqref{resolution of subregular}, we obtain a resolution of $T^{\ell,-}_\lambda$:
\begin{align*}
0\rightarrow T_\lambda^{\ell,-}\rightarrow F^k_{\lambda,0}\rightarrow F^k_{\lambda,1}\rightarrow \cdots \rightarrow 0.
\end{align*}
with 
$$F^k_{\lambda,p}\simeq \bigoplus_{l(w)=p}S_{(w^{-1}*\lambda-\lambda,\varpi_{n-1})} W_{w^{-1}*\lambda}^\ell.$$
Here $S_{\theta}$ ($\theta\in \Z$) are the spectral flow twists given\cite{Li} by 
\begin{align*}
&S_{\theta}\colon X(z)\mapsto Y(\Delta(\theta H_-,z)X,z),\quad (X\in\W^\ell(\sll_{n|1},f_{\mathrm{prin}})),\\
&\Delta(\theta H_-,z):=z^{\theta H_{-,0}}\mathrm{exp}\left(-\theta \sum_{m=1}^\infty \tfrac{H_{-,m}}{m}(-z)^{-m}\right).
\end{align*}
The character $\mathrm{tr}_\bullet q^{L_0}z^{H_{-,0}}$ of $W_\lambda^\ell$ is given by
\begin{align*}
\mathrm{ch}\ W^{\ell,-}_\lambda=q^{\Delta_\lambda^k+\frac{n}{2}(\lambda,\varpi_{n-1})-\frac{1}{2}\epsilon^2(\lambda,\varpi_{n-1})^2}z^{\epsilon^2(\lambda,\varpi_{n-1})}
\frac{\poch{-zq^{\frac{n+1}{2}},-z^{-1}q^{\frac{-n+3}{2}}}}{\poch{q}^n}.
\end{align*}
Let us compare it with the characters of the Wakimoto modules over $\W^\ell(\sll_{n|1},f_{\mathrm{prin}})$. These modules are defined as 
$$\mathbb{W}^{\ell,-}_{\mu}=V_{\Z}\otimes\pi^{\ell+n-1}_{\h,\mu},\quad (\mu\in\h^*),$$
through the free field realization \eqref{FFR}. Since \eqref{FFR} sends
\begin{align*}
&L(z)\mapsto (:\partial c(z)\ b(z):+\tfrac{n-1}{2}\partial :c(z)b(z) :)\otimes 1\\
&\hspace{2cm}+1 \otimes (L_{\mathrm{sug}}(z)+(1-\tfrac{1}{\ell+n-1})\partial \rho(z)+\tfrac{1-n}{2}\partial \varpi_n(z)),\\
&H_-(z)\mapsto :b(z)c(z):\otimes 1+1\otimes \varpi_n(z),
\end{align*} 
their characters $\mathrm{tr}_\bullet(q^{L_0}z^{H_-,0})$ are given by 
\begin{align}
\mathrm{ch}\ \mathbb{W}^{\ell,-}_{\mu}
=q^{\Delta_{\mu}^\ell+\frac{n-1}{2}(\mu,\varpi_n)}z^{(\mu,\varpi_n)}\frac{\poch{-zq^{\frac{n+1}{2}},-z^{-1}q^{\frac{-n+1}{2}} }}{\poch{q}^n }
\end{align}
with 
$\Delta_\mu^\ell=\tfrac{1}{2(\ell+n-1)}(\mu,\mu)-(1-\tfrac{1}{\ell+n-1})(\mu,\rho).$
Then we find that 
\begin{align*}
\mathrm{ch}\ W^\ell_\lambda=\frac{1}{1+z^{-1}q^{\frac{-n+1}{2}} }\mathrm{ch}\ \mathbb{W}^{\ell,-}_{\mu},\quad \mu=-(\ell+n-1)(\tilde{\lambda}+(\lambda,\varpi_{n-1})\varpi_n).
\end{align*}
Here $\tilde{\lambda}$ is the image of $\lambda$ under the orthogonal decomposition 
$\h(\sll_{n|1})=\h(\sll_n)\oplus \C \varpi_n.$ 
Indeed, $W^\ell_\lambda$ is a submodule of $\mathbb{W}^{\ell,-}_\mu$ (``thin Wakimoto module''). This is proven by introducing ``thick'' Wakimoto modules 
$$\widehat{\mathbb{W}}^{k,+}_\lambda:=\pi^{k+n}_{\h,\lambda}\otimes \Pi,\quad (\lambda\in \h^*),$$
whose image coincide with the Wakimoto modules $\mathbb{W}^{\ell,-}_\mu$. 
The $\W^k(\sll_n,f_{\mathrm{sub}})$-module structure on $\widehat{\mathbb{W}}^{k,+}_\lambda$ is given by the composition of \eqref{Wakimoto realization of subregular} and the bosonization $\beta\gamma\rightarrow \Pi$, that is, the free field realization \eqref{FFR}. Replacing $\Pi$ with its modules $\Pi_{0}[a]$ ($a\in \C$) below, we obtain the general Wakimoto module for $\W^\ell(\sll_{n|1},f_{\mathrm{prin}})$.

To obtain Verma-type modules on $\W^\ell(\sll_{n|1},f_{\mathrm{prin}})$-side as in the case of $V^c(\ns_2)$, we need to replace the Fock modules $\pi^{k+n}_{\h,\lambda}$ by Verma modules for $\W^k(\sll_n)$ and thus need the twist \eqref{FFR} and use \eqref{invHR in general}, which satisfies
\begin{align*}
\begin{array}{cll}
\W^k(\sll_n,f_{\mathrm{sub}})&\rightarrow&\hspace{1cm}\W^k(\sll_n)\otimes \Pi\\
L(z)&\mapsto&W_2(z)\otimes1 +1\otimes (\tfrac{1}{2}:c(z)d(z) :-\tfrac{n}{2}\partial \mathbf{b}^-(z))\\
H_-(z)&\mapsto &\hspace{1cm} 1\otimes \mathbf{b}^+(z)
\end{array}
\end{align*}
where $W_2(z)$ is the Virasoro field of $\W^k(\sll_n)$ and $\mathbf{b}^\pm(z)$ are given by
$$\mathbf{b}^\pm(z)=\pm \tfrac{1}{2}(\varepsilon_+-1)c(z)+\tfrac{1}{2}d(z).$$
Let $\mathbb{M}^k_{\chi_\lambda}$ ($\lambda\in \h^*$) be the Verma module of $\W^k(\sll_n)$ equipped with the central character $\chi_\lambda\colon \mathcal{Z}(\sll_n)\rightarrow \C$. It is generated by the one dimensional representation $\chi_\lambda$ of the Zhu's algebra $\mathrm{Zhu}(\W^k(\sll_n))\simeq \mathcal{Z}(\sll_n)$. 
Under this notation, we have a natural homomorphism $\mathbb{M}^k_{\chi_\lambda}\rightarrow \pi^{k+n}_{\h,\lambda}$ of $\W^k(\sll_n)$-modules, which sends the highest weight vector to the other. Let $\Pi_\theta[a]$ $(\theta,a\in \C)$ be the $\Pi$-module defined by 
$$\Pi_\theta[a]=\bigoplus_{m\in \Z}\pi^{\Z+\ssqrt{-1}\Z}_{\theta\mathbf{b}^++(m+a)c}.$$
Then we define the relaxed highest weight module $R_\theta^k(\chi_\lambda,a)$ of $\W^k(\sll_n,f_{\mathrm{sub}})$ as
$$R_\theta^{k,+}(\chi_\lambda,a):=\mathbb{M}^k_{\chi_\lambda}\otimes \Pi_\theta[a].$$
It is positively graded if and only if $\theta=-n/2$ and we set $R^{k,+}(\chi_\lambda,a)=R^{k,+}_{-n/2}(\chi_\lambda,a)$, whose character is given by
\begin{align*}
\mathrm{ch}\ R^{k,+}(\chi_\lambda,a)=q^{\Delta^k_\lambda+\frac{1}{8}n^2\varepsilon_+}z^{a-\frac{n}{2}\varepsilon_+}\frac{\sum_{m\in\Z}z^m}{\poch{q}^{n+1}}.
\end{align*} 
Applying the functor $H^{\mathrm{rel}}_{+,\epsilon(a-\frac{n}{2}\varepsilon_+)}$, we find that the image has the following character 
:\begin{align*}
&\mathrm{ch}\ H^{\mathrm{rel}}_{+,\epsilon(a-\frac{n}{2}\varepsilon_+)}(R^{k,+}(\chi_\lambda,a))\\
&\hspace{1cm}=q^{\Delta_\lambda^k+\frac{1}{8}n^2 \varepsilon_+^2-\frac{1}{2}\epsilon^2(a-\frac{n}{2}\varepsilon_+)^2} z^{\epsilon^2(a-\frac{n}{2}\varepsilon_+)}\frac{\poch{-zq^{\frac{1}{2}},-z^{-1}q^{\frac{1}{2}}} }{\poch{q}^n }.
\end{align*} 
When $n$ is even, the image is given by the spectral flow twist $S_{-n/2}(R^{\ell,-}_{\chi_{\mu}})$ with $\mu=-(\ell+n-1)(\tilde{\lambda}+a\varpi_n)$ of the Verma-type $\W^\ell(\sll_{n|1},f_{\mathrm{prin}})$-module $R^{\ell,-}(\chi_{\mu})$, which is defined in general as
$$R^{\ell,-}(\chi_\mu)=\mathbb{M}^\ell_{\chi_\mu}\otimes V_\Z,\quad (\mu\in \h^*),$$
where $\mathbb{M}^\ell_{\chi_\mu}$ is the Verma module of $\W^\ell(\gl_n)$ with a natural homomorphism $\mathbb{M}^\ell_{\chi_\mu}\rightarrow \pi^{\ell+n-1}_{\h,\mu}$ and $\W^\ell(\sll_{n|1},f_{\mathrm{prin}})$ acts on $R^{\ell,-}(\chi_\mu,a)$  through $\W^\ell(\sll_{n|1},f_{\mathrm{prin}})\rightarrow \W^\ell(\gl_n)\otimes V_\Z$.
The spectral flow twist makes the above asymmetric Verma-type module into a ``(symmetric) Verma module''.
To summarize, we have obtained the following dictionary generalizing Tab.\ \ref{correspondence of modules in general} below.
\renewcommand{\arraystretch}{1.3}
\begin{table}[h]
\caption{}
\label{correspondence of modules in general}
\begin{tabular}{cc}
\hline
$\W^k(\sll_n,f_{\mathrm{sub}})$-side & $\W^\ell(\sll_{n|1},f_{\mathrm{prin}})$-side \\ \hline
Wakimoto modules & ``thin'' Wakimoto modules  \\
``thick'' Wakimoto modules & Wakimoto modules  \\
relaxed highest weight modules & ``Verma modules'' \\ \hline
\end{tabular}
\end{table}
\renewcommand{\arraystretch}{1}

It is an interesting problem to extend the correspondence of classes of modules in the spirit of the quantum Langlands duality\cite{AF} of Arakawa and Frenkel for modules over the princial $\W$-algebras.
On the other hand, since the positively graded modules are controlled by the finite $\W$-superalgebras\cite{DK}, it is natural to compare the finite $\W$-superalgebras in our case. In the literature\cite{EGG,LT}, the finite $\W$-algebras for $\W^k(\sll_{n}, f_{\mathrm{sub}})$ have a nice realization as infinitesimal Cherednik algebras.

\section{Conclusions}

The pair of $\W$-superalgebras $\W^k(\sll_n,f_{\mathrm{sub}})$ and $\W^\ell(\sll_{n|1},f_{\mathrm{prin}})$ enjoy the Feigin--Frenkel type duality for their Heisenberg cosets. It can be proven either by using the uniqueness property of the universal $\W_\infty$-algebra $\W_\infty[c,\lambda]$ or by the technique of free field realization. Furthermore, the duality is refined to a reconstruction type theorem of one of the $\W$-superalgebras from the other by the coset functor or the relative semi-infinite cohomology functor.
Either of them leads to the block-wise equivalence of module categories consisting of Kazhdan--Lusztig objects with respect to the Heisenberg subalgebra. This equivalence is incorporated with isomorphisms of all the logarithmic intertwining operators. Moreover, we have obtained a dictionary of the correspondence of some fundamental class of modules.
There are several interesting problems still left open.
\begin{itemize}
\item The correspondence for correlation functions on the Riemannian sphere has been established\cite{CHS} by path integral method, which is worth generalization to correlation functions on arbitrary Riemannian surfaces.
\item Apart from the rational cases, $\W_k(\sll_n,f_{\mathrm{sub}})$ at $k=-n+\frac{n}{n+1}$, known also as the chiral algebras of Argyres--Douglas theories of type $(A_1,A_{2n-1})$\cite{ACGY} have suitable module categories which have a non-semisimple braided tensor category structure \cite{ACKR}. It is an interesting problem to understand the $\W_\ell(\sll_{n|1},f_{\mathrm{prin}})$-side and compare the braided tensor category structure.
\end{itemize}

Obviously, it is interesting to pursue the correspondence of modules and intertwining operators, or correlation functions in general, for other pairs of $\W$-superalgebras, namely hook-type $\W$-superalgebras appearing in Refs. \citen{GR,CL1,CL2}. 
The Feigin--Frenkel type duality is already established\cite{CL1,CL2} mathematically, the reconstruction type theorem is also established for a large class\cite{CLNS} and the correspondence of module categories and the behavior of $q$-characters are work in progress.\cite{CLN}

\bibliographystyle{plain}
\bibliography{nakatsuka}
\end{document}